\documentclass[12pt,draft]{article}
\usepackage{amsmath,amsfonts,amsthm,amssymb,amscd}
\newcommand{\llangle}{\langle\langle}
\newcommand{\rrangle}{\rangle\rangle}

\newcommand{\p}{\partial}

\newcommand{\e}{\varepsilon}
\newcommand{\vk}{\varkappa}
\newcommand{\tauu}{t}
\newcommand{\vp}{\varphi}
\newcommand{\om}{W}
\newcommand{\m}{{\bf m}}
\newcommand{\vv}{{\bf v}}
\newcommand{\vvv}{{\tilde{\bf v}_k^\gamma}}
\newcommand{\R}{{\mathbb R}}
\newcommand{\C}{{\mathbb C}}
\newcommand{\IP}{{\bf P}}
\newcommand{\Z}{{\mathbb Z}}
\newcommand{\E}{{\bf E}}
\newcommand{\cE}{{\cal E}}
\newcommand{\T}{{\mathbb T}}
\newcommand{\N}{{\mathbb N}}
\newcommand{\PP}{{\bf P}}
\newcommand{\II}{{I}}

\newcommand{\fA}{{\frak A}_{\alpha_j}}
\newcommand{\BB}{{\cal B}}
\newcommand{\DD}{{\cal D}}
\newcommand{\FF}{{\cal F}}
\newcommand{\HH}{{h}}
\newcommand{\QQ}{{\cal Q}}

\newcommand{\strela}{\rightharpoonup}

\newcommand{\Plip}{{\mathop{\rm Lip}\nolimits}{}_{\mathop{\rm
Lock}\nolimits}}
\newcommand{\const}{\mathop{\rm const}\nolimits}

\newcommand{\diver}{\mathop{\rm div}\nolimits}

\def\S{\mathhexbox278}
\def\lan{\langle}
\def\ran{\rangle}

\theoremstyle{plain}
\newtheorem{theorem}{Theorem}[section]
\newtheorem{lemma}[theorem]{Lemma}
\newtheorem{proposition}[theorem]{Proposition}
\newtheorem{corollary}[theorem]{Corollary}
\theoremstyle{definition}

\theoremstyle{remark}

\numberwithin{equation}{section}

\begin{document}
\author{Sergei B. Kuksin and Andrey L. Piatnitski}
\title{Khasminskii--Whitham averaging for randomly
perturbed KdV equation.}
\date{}
\date{(to appear in JMPA)}
\maketitle

\begin{abstract}
We consider the damped-driven KdV equation
$$
\dot u-\nu{u_{xx}}+u_{xxx}-6uu_x=\sqrt\nu\,\eta(t,x),\; x\in S^1,\
\int u\,dx\equiv \int\eta\,dx\equiv0\,,
$$
where $0<\nu\le1$ and the random process $\eta$ is smooth in $x$ and
white in $t$. For any periodic function $u(x)$ let $
I=(I_1,I_2,\dots) $ be the vector, formed by the KdV integrals of
motion,
 calculated for the potential $u(x)$. We prove that
 if $u(t,x)$ is a solution of the equation above,
 then for $0\le t\lesssim\nu^{-1}$ and $\nu\to0$ the vector
$
I(t)=(I_1(u(t,\cdot)),I_2(u(t,\cdot)),\dots)
$
satisfies the (Whitham) averaged equation.
\end{abstract}

\tableofcontents

\setcounter{section}{-1}
\section{Introduction}\label{s0} It is well known since the pioneer
works of Novikov and Lax that the KdV equation
\begin{equation}\label{0.1}
    \dot u+u_{xxx}-6uu_x=0,  
\end{equation}
defines an integrable infinite dimensional Hamiltonian system in a
space $H^p$ of $2\pi$-periodic Sobolev functions of order $p\ge0$
with zero meanvalue. It means that  KdV  has infinitely many
integrals of motion $I_1,I_2,\dots$, which are non-negative analytic
functions on $H^p$, and for any non-negative sequence
$I=(I_1,I_2,\dots)$
 the set
$T_I=\{u:\, I_j(u)=I_j\ \forall\,j\}$ is an analytic torus in $H^p$
of dimension $|J(I)|\le\infty$, where $J$ is the set $J=\{j:\,
I_j>0\}$. Each torus $T_I$ carries an analytic cyclic coordinate
$\vp=\{\vp_j,j\in J(I)\}$, and in the coordinates $(I,\vp)$ the
KdV-dynamics takes the integrable form
\begin{equation}\label{0.2}
    \dot I=0,\quad \dot\vp=W(I)\,.
\end{equation}
The {\it frequency vector } $W$ analytically depends on $I$. See
\cite{McKT76,KaP} and Section~\ref{s2} below.

Importance of these remarkable features of  KdV  is jeopardised by
the fact that KdV arises in physics only as an approximation for
`real' equations, and it is still unclear up to what extend the
integrability property persists in the `real' equations, or how it
can be used to study them.

The persistence problem turned out to be difficult, and the progress
in its study is slow. In particular, it was established that small
Hamiltonian perturbations of  KdV  do not destroy
majority of time-quasiperiodic solutions, corresponding to
\eqref{0.2} with $|J(I)|<\infty$ (see \cite{K2,KaP}), but it is
unknown how these perturbations affect the almost-periodic solutions
($|J(I)|=\infty$), and whether solutions of the perturbed equations
are stable in the sense of Nekhoroshev.

Probably it is even more important to understand the behaviour of
solutions for KdV, perturbed by non-Hamiltonian terms
(e.g., to understand how small dissipation affects the equation).
The first step here should be to study how a $\nu$-perturbation
 affects the dynamics \eqref{0.2} on time-intervals of
order $\nu^{-1}$. For perturbations of finite-dimensional integrable
systems this question is addressed by the classical averaging
theory, originated by Laplace and Lagrange. During more than 200
years of its history this theory was much developed, and good
understanding of the involved phenomena was achieves, e.g. see in
\cite{AKN}. In particular, it is known that for a perturbed
finite-dimensional integrable system
\begin{equation}\label{0.3}
    \dot I=\nu f(I,\vp)\quad \dot\vp =W(I)+\nu g(I,\vp),
    \quad\nu\ll1,
\end{equation}
where $I\in\R^n,\ \vp\in\T^n$, 
 on time-intervals of order $\nu^{-1}$ the action $I(t)$
 may be well approximated by solutions of the {\it averaged
 equation}
 \begin{equation}\label{0.4}
    \dot I=\nu\lan f\ran(I),\quad \lan
    f\ran(I)=\int_{\T^n}f(I,\vp)\,d\vp\,,
 \end{equation}
 provided that the initial data
$(I(0),\vp(0))$ are typical. This assertion is known as the {\it
averaging principle}.

The behaviour of solutions of infinite-dimensional systems \eqref{0.3}
on time-inter\-vals of order $\gtrsim\nu^{-1}$ is poorly  understood.
 Still applied mathematicians believe that the averaging
principle holds, and use \eqref{0.4} to study solutions of
\eqref{0.3} with $n=\infty$. In particular, if \eqref{0.3} is a
perturbed KdV equation, written in the variables $(I,\vp)$, then
\eqref{0.4} is often called the {\it Whitham equation}
(corresponding to the perturbed KdV). The approximation for $I(t)$
in \eqref{0.3} with $0\le t\le\nu^{-1}$ by $I(t)$, satisfying
\eqref{0.4}, is called the {\it Whitham averaging principle} since
in \cite{Wh} the averaging is systematically used in similar
situations. In so far the Whitham averaging for the perturbed KdV
equation under periodic boundary conditions was not rigorously
justified. Instead mathematicians, working in this field, either
postulate the averaging principle and study the averaged equations
(e.g., see \cite{FFM} and \cite{DN89}), or postulate that the
solution regularly -- in certain sense -- depends on the small
parameter and show that this assumption implies the Whitham
principle, see \cite{Kri88}.

 The main goal of this paper is to justify the Whitham averaging
  for randomly perturbed  equations.

Let us start with random perturbations of the integrable system
\eqref{0.2} with $I\in\R^n, \vp\in\T^n$, where $n<\infty$.
Introducing the fast time $\tau=\nu t$ we write the perturbed
system as the Ito equation
\begin{equation}\label{0.5}
\begin{split}
dI&=F\,d\tau+\sigma\,d\beta_\tau,    \\
d\vp&=(\nu^{-1}W(I)+G)\,d\tau+g\,d\beta_\tau\,.
\end{split}
\end{equation}
Here $F,G,\sigma$ and $g$ depend on $(I,\vp)$, $\beta_\tau$ is a
vector-valued Brownian motion and $\sigma, g$ are matrices. It was
claimed in \cite{Khas68} \footnote{The main theorem of \cite{Khas68}
deals with the situation when the unperturbed system is a stochastic
equation with a non-degenerate diffusion for $\vp$, but in its last
section it is claimed that the ideas of the proof also apply to
\eqref{0.5}.} and proved in \cite{WF03} that (under certain
assumptions, where the main one is non-degeneracy of the diffusion
$\sigma$ and of the frequency-map $W$)  when $\nu\to0$, the solution
$I(\tau)$ converges in distribution to a solution of the averaged
equation
\begin{equation}\label{0.6}
    dI=\lan F\ran(I)\,d\tau+\lan\sigma\ran(I)\,d\beta\,,
\end{equation}
where $\lan F\ran$ is defined as in \eqref{0.4} and the matrix
$\lan\sigma\ran(I)$ is a symmetric square root of the matrix
$\int_{T^n}\sigma\sigma^t\,d\vp$.\smallskip

Now let us consider a randomly perturbed (`damped-driven') KdV
equation
\begin{equation}\label{1}
\dot u-\nu{u_{xx}}+u_{xxx}-6uu_x=\sqrt\nu\eta(t,x)\,.
\end{equation}
As before, $x\in S^1$ and  $\int u\,dx\equiv\int\eta\,dx\equiv0$.
The force  $\eta$ is a Gaussian random field, white in time $t$:
$$
\eta=\frac{\partial}{\partial
t}\sum_{s\in\Z_0}b_s\beta_s(t)e_s(x)\,,
$$
where $\Z_0=\Z\setminus\{0\}$, $\beta_s(t)$ are standard independent
Wiener processes, and $\{e_s, s\in\Z_0\}$ is the usual trigonometric
basis
\begin{equation}\label{3}
 e_s(x)=
 \begin{cases}
   \cos s x, \quad s>0\,,\\
  \sin s  x, \quad s <0\,.
  \end{cases}
\end{equation}
Concerning the real constants $b_s$ we assume that
\begin{equation}\label{4}
    b_s\le C_m|s|^{-m}\quad \forall\,m,s
\end{equation}
with some constants $C_m$ (so $\eta(t,x)$ is smooth in $x$), and
\begin{equation}\label{5}
    b_s\ne0\quad\forall\,s\,.
\end{equation}
The factor $\sqrt\nu$ in front of the force $\eta (t,\, x)$ is
natural since under this scaling solutions of \eqref{1} remains of
order 1 as $t \to \infty$ and $\nu \to 0$. Eq. \eqref{1} defines a
Markov process in the function space $H^p$. Due to \eqref{5} it has
a unique stationary measure. Let $u^\nu (t,\, x )$, $t\geq 0$, be a
corresponding stationary in time solution for \eqref{1}; or let
$u^\nu$ be a solution, satisfying
\begin{equation}\label{0.7}
u^\nu (0,\, x) = u_0 (x),
\end{equation}
where $u_0(x)$ is a non-random smooth function. In Section~\ref{s1}
we prove that all moments of all Sobolev norms $\| u^\nu (t,\,
\cdot)\|_m$ are bounded uniformly in $\nu >0$ and $t\geq 0$. Let us
write $u^\nu (\tau)$ as $(I^\nu (\tau),\, \vp^\nu (\tau ))$. These
processes satisfy the infinite-dimensional equation \eqref{0.5}, so
by the just mentioned estimates the processes $\{I^\nu
(\cdot),0<\nu\le1\}$ form a tight family, and along suitable
sequences $\nu_j \to 0$ we have a weak convergence in distribution
\begin{equation}\label{0.8}
I^{\nu_j} (\cdot) \to I^0(\cdot),
\end{equation}
where, according to the type of the solutions $u^\nu (\tau)$, the
limiting process $I^0 (\tau)$ is either stationary in $\tau$, or
satisfies $I^0 (0) = I(u_0(\cdot))$.

The main results of this work are the following two theorems, proved
in Section \ref{s6}:\smallskip

\noindent{\bf Theorem A.}
 The limiting process $I^0(\tau)$ satisfies the Whitham
equation \eqref{0.6}, corresponding to the perturbed KdV equation
\eqref1. It is non-dege\-ne\-rate in the sense that for any $\tau>0$
and each $k\ge1$ we have $\PP\{I_k^0(\tau)=0\}=0$.
\smallskip

\noindent{\bf Theorem B.} If the processes $u^\nu (\tau)$ are
stationary in $\tau$, then for any $\tau \geq 0$ the law of the pair
$(I^{\nu_j}(\tau),\, \vp^{\nu_j} (\tau))$ converges to the product
measure $q^0\times d\vp$, where $q^0$ is the law of $I^0(0)$ and
$d\vp$ is the Haar measure on $\T^\infty$.\smallskip

The proof is based on the scheme, suggested by Khasminskii in
\cite{Khas68}, see also \cite{FW98} and \cite{Ver91}.
  It uses the estimates from
Section~\ref{s1} and more sophisticated estimates, obtained in
Sections~\ref{s5} and \ref{s.av}.  Namely, we use crucially
Lemma~\ref{notatzero} (Section~\ref{s5}) and  Lemma \ref{l.new}
(Section~\ref{s.av}). In the former  coupling arguments are evoked
to prove that  for any $k$ probability of the event $\{I_k^\nu (t)
<\delta\}$ goes to zero with $\delta$,   uniformly in $\nu$ and $t$.
This is important since \eqref{0.5} is an equation for $I$ in the
octant $\{ I\, | \,I_j
> 0\quad \forall j\}$ which degenerates at the boundary $\{ I \,| \,I_j =
0 \ \text{for some}\ j\}$. In the latter  we
 examine the random process $W^m (\tau) = W^m
(I^\nu (\tau))$, where $W^m$ is the vector, formed by the first $m$
components of the frequency vector $W$. Exploiting Krylov's results
from \cite{Kr77} we estimate the density against the Lebesgue
measure of the law of the averaged vector $s^{-1} \int_0^s W^m
(I^\nu(\tau)) \, d\tau$, $s\sim 1$. We use this estimate to show
that with
 probability close to one the components of the  vector $W^m (\tau)$ are
non-commensurable, so the fast motion $(d/d\tau)\vp^m = \nu^{-1}
W^m(\tau)$ is ergodic on the torus $\T^m \subset \T^\infty$, for any
$m$. This is a crucial step of the proof of Theorem~A. Our proof of
Lemma \ref{l.new} is `hard' in the sense that it uses heavily the
analyticity of the frequency map $W(I)$.

The arguments above are applied to the perturbed KdV equation,
written in the Birkhoff normal form (eq.~\eqref{2.2}  in
Section~\ref{s2}). They apply as well to perturbations of other
Birkhoff-integrable
 equations if their solutions satisfy good apriori
 estimates uniformly in the small parameter, and the corresponding
 transformation to
 the Birkhoff coordinates is smooth and is polynomially bounded at
infinity. In the KdV case which we consider, half of the required
 bounds on the transformation
   is established in the recent paper \cite{Kor06}. We are
certain that the remaining half can be obtained similarly, but do
 not prove them in this work, see Theorem~\ref{t2.2} in Section~\ref{s2}.
\smallskip

The Whitham equation \eqref{0.6}, corresponding to the perturbed KdV
\eqref1, is a complicated infinite-dimensional stochastic
differential equation.
 Theorem~A implies that for any smooth initial
data $I(0)$ it has a  weak solution, but we do not know if
this solution is unique.
 We point out that, firstly, if \eqref{0.6}
has a unique solution and the process $u^\nu (\tau)$ satisfy
\eqref{0.7}, then the law of the limiting process $I^0$ is
independent of the sequence $\{\nu_j\}$, and the convergence
\eqref{0.8} holds for $\nu\to 0$. Secondly, if \eqref{0.6} has a
unique stationary measure, then a similar assertion holds for
stationary solutions $u^\nu (\tau)$.\medskip

\textbf{The inviscid limit}. Let us consider the stationary
solutions of eq. \eqref{1} in the original time $t$. The apriori
estimates from Section~\ref{s1} imply that this family is tight
in $C([0,T];H^p)$ for any $T>0$ and any $p>0$.
Therefore, along sequences $\nu_j \to 0$, we have convergence in
distribution
\begin{equation}\label{0.10}
u^{\nu_j} (\cdot) \to u^0 (\cdot)
\end{equation}
(the limiting process $u^0(t)$ apriori depends of the sequence
$\{\nu_j\}$). The arguments, applied in Section~10 of \cite{K3} to
the randomly perturbed Navier - Stokes equation \eqref{NSE} also
apply to \eqref1. They show that a.e. realisation of the limiting
process $u^0 (t,\, x)$ is a smooth solution of the KdV equation
\eqref{0.1}. In particular, the law $\mu^0$ of the random variable
$u^0(0,\cdot)\in H^p$ is an invariant measure for the dynamical
system which  KdV defines in $H^p$.  But KdV has infinitely many
integrals of motion; so it has a lot of invariant measures. How to
distinguish among them the measure $\mu^0$? Noting that
$u^\nu(t)_{t=0} = u^{\nu}(\tau)_{\tau =0}$, we apply Theorem~B to
get that the isomorphism $u(\cdot)\mapsto (I,\, \vp )$ transforms
$\mu^0$ to the measure $q^0\times d\vp$. In particular, if
\eqref{0.6} has a unique stationary measure, then the measure
$\mu^0$ is uniquely defined, and the convergence \eqref{0.10} holds
for $\nu\to 0$.

This discussion shows that in difference with the deterministic
situation, averaged randomly perturbed equations describe not only
behaviour of solutions for a pre-limiting equation on time-intervals
of order $\nu^{-1}$, but also its asymptotic in time properties.
Indeed, under the double limit `first $t\to\infty$,   next
$\nu\to0$', the distribution of any solution converges to a measure,
simply expressed in terms of a stationary measure of the averaged
equation.
\medskip

\textbf{The Eulerian limit}. The perturbed KdV equation \eqref{1} is
a reasonable model for the randomly perturbed 2D NSE
\begin{equation}\label{NSE}
\begin{split}
\dot u -\nu\Delta u +(u\cdot\nabla)u + \nabla p = \sqrt\nu \eta
(t,\,
x),\,\, x\in\T^2,\\
\diver u =0,\,\, \int u\, dx \equiv \int\eta \, dx\equiv 0,
\end{split}
\end{equation}
obtained by replacing in \eqref{NSE}
 the 2D Euler equation \eqref{NSE}${}_{\nu=0}$ (which is a
Hamiltonian PDE with infinitely many integrals of motion) by KdV.
Under restrictions on the random force  $\eta (t,\,x)$, similar to
those imposed on the force in \eqref{1}, eq.~\eqref{NSE}
(interpreted as a Markov process in the space of divergence-free
vector fields $u(x)$), has a unique stationary measure, see in
\cite{K3}. Let $(u^\nu (t),\, p^\nu (t))$ be the corresponding
stationary solution. Then, along sequences $\nu_j\to 0$, the
convergence in distribution holds
\begin{equation}
(u^{\nu_j}(\cdot),\, p^{\nu_j}(\cdot))\to (u^0 (\cdot ),\, p^0
(\cdot )),
\end{equation}
where the limiting process $(u^0,\, p^0)$ is stationary in time, is
sufficiently smooth in $t$ and $x$, and a.e. its realisation
satisfies the free Euler equation \eqref{NSE}${}_{\nu=0}$.
Accordingly, the law $\mu^0$ of $u^0(0)$ is an invariant measure for
the dynamical system, which the Euler equation defines in the space
of divergence-free vector fields. To study the measure $\mu^0$ (in
fact, the set of measures $\mu^0$, since it is possible that now the
limit depends on the sequence $\{\nu_j\}$), is an important problem
in (mathematical) $2D$ turbulence. The problem, addressed in this
work, may be considered as its model.

\bigskip
\noindent {\it Agreements.} Analyticity of maps $B_1\to B_2$ between
Banach spaces $B_1$ and $B_2$, which are the real parts of complex
spaces $ B_1^c$ and $B_2^c$, is understood in the sense of
Fr\'echet. All analytic maps which we consider possess the following
additional property: for any $R$ a map analytically extends to a
complex $(\delta_R>0)$--neighbourhood of the ball $\{|u|_{B_1}<R\}$
in $B_1^c$.   When two random variables are equal
 almost sure, we usually drop the specification ``a.s.".

\noindent {\it  Notations.}  $\chi_A$ stands for the indicator
function of a set $A$ (equal 1 in $A$ and equal 0 outside $A$). By
$\vk(t)$ we denote various functions of $t$ such that $\vk(t)\to0$
when $t\to\infty$, and by $\vk_\infty(t)$ denote functions $\vk(t)$
such that $\vk(t)=o(t^{-N})$ for each $N$. We write
$\vk(t)=\vk(t;R)$ to indicate that $\vk(t)$ depends on a parameter
$R$. For a measurable set $Q\subset\R^n$ we denote by $|Q|$ its
Lebesgue measure.

\section{The equation and its solutions}\label{s1}

We denote by $H$ the Hilbert space
$$
H=\{u\in L_2(S^1):\, \int u\,dx=0\}
$$
with the scalar product $\  \langle
u,v\rangle=\frac1{\pi}\int_0^{2\pi}u(x)v(x)\,dx. $ Then $\{e_s, s\in
\Z_0\}$ (see \eqref{3}) is its Hilbert basis. We set $H^m$ to be the
$m$-th Sobolev space, formed by functions with zero mean--value, and
given the norm $\|u\|_m=\langle \frac{\p^mu}{\p x^m},
\frac{\p^mu}{\p x^m}\rangle^{1/2}$.

We write  the KdV equation as
\begin{equation}\label{2.1}
\dot u+V(u)=0\,,\quad V(u)=u_{xxx}-6uu_x\,,
\end{equation}
and  re-write eq. \eqref{1} as
\begin{equation}\label{6}
    \dot u-\nu {u_{xx}}+V(u)=\sqrt\nu\,\eta(t,x)\,.
\end{equation}

It is well known that a dissipative nonlinear equation in one
space--dimension with a white in time r.h.s. has a unique strong
solution if the equation's solutions satisfy sufficiently strong
a-priori estimates. In Appendix we show that any smooth solution of
\eqref{1} with a deterministic initial data
\begin{equation}\label{1.14}
u(0)=u_0\,,
\end{equation}
where $u_0\in H^m, m\ge1,$ satisfies the following estimates:
\begin{equation}\label{1.4}
\E e^{\sigma \| u (t)\|^2_0}\le \max \big(\E e ^{\sigma \|u
(0)\|^2_0},\,\, 2e^{2\sigma B_0}\big),
\end{equation}
\begin{equation}\label{1.12}
\E\|u(t)\|_m^2\le\max\big(4\E\|u(0)\|^2_m,C'_m\big),
\end{equation}
\begin{equation}\label{1.13}
\E\|u(t)\|_m^k\le C\big(\|u_0\|_{m k},B_{m+1},m,k\big)\,.
\end{equation}
Here $t\ge0$, $k\in\N$ and $\sigma\le(2\max{b_s^2})^{-1}$.

Accordingly, we have the following result:
\begin{theorem}\label{t1.1} For any deterministic $u_0\in H^m$,
$m\ge1$, the problem \eqref{1}, \eqref{1.14} has a unique solution
$u(t,x)$. It satisfies estimates \eqref{1.4} -
\eqref{1.13}.
\end{theorem}

Due to assumption \eqref{5}, eq.~\eqref{1} has a unique stationary
measure $\mu_\nu$ and any solution converges to $\mu_\nu$ in
distribution. For  the randomly forced 2D~NSE equation this result
now is well known (e.g., see in \cite{K3}). The proofs for
eq.~\eqref{1} are  simpler and we do not discuss them.

Let $u^0_\nu(t,x)$ be a solution of \eqref{1}, \eqref{1.14} with
$u_0=0$. Since $\DD(u^0_\nu(t))\strela\mu_\nu$, then
Theorem~\ref{t1.1} and the Fatou lemma imply

\begin{theorem}\label{tI.1}
The unique stationary measure $\mu_\nu$ satisfies the estimates
\begin{equation*}\begin{split}
    &\int_He^{\sigma\|u\|_0^2}\mu_\nu(du)\le C_\sigma<\infty
\quad \forall\,\sigma\le(2\max{b_s^2})^{-1}
,\\
    &\int_H\|u\|_m^k \mu_\nu(du)\le
    C_{m,k}<\infty\quad\forall\,m\,,k\,.
\end{split}\end{equation*}
\end{theorem}

\section{Preliminaries on the KdV equation}\label{s2}
In this section we discuss integrability of the KdV equation
\eqref{2.1}.

 For $r\ge0$ let us denote by $\HH^r$ an abstract Hilbert space
with the basis $\{f_j,j=\pm1,\pm2,\dots\}$ and the norm $|\cdot|_r$,
where
$$
|v|_r^2=\sum_{j\ge1}j^{1+2r}(v_j^2+v_{-j}^2)\quad\text{for}\quad
v=\sum_{j\in\Z_0}v_jf_j.
$$
We denote  $\vv_j=  \left(\begin{array}{c}
    v_j\\
    v_{-j}\\
  \end{array}\right)$,
  and identify a
vector $v=\sum v_jf_j\in\HH^r$ with the sequence
$(\vv_1,\vv_2,\dots)$.

\begin{theorem}\label{t2.1} {\rm (see \cite{KaP})}.  There exist an
analytic diffeomorphism $\Psi:H\to\HH^0$ and an analytic functional
$K$ on $\HH^0$ of the form
$$
K(\sum v_jf_j)=\tilde
K(I_1,I_2,\dots)\,,\;\;I_j=\frac12\,(v_j^2+v_{-j}^2)\,,
$$
with the following properties:

1)  $\Psi$ defines, for any $m\in\N$,  an analytic diffeomorphism
$\Psi:H^m\to\HH^m$,

2)  $d\Psi(0)$ is the map $\  H^m\ni\sum u_s e_s\mapsto \sum
|s|^{-1/2}v_sf_s\in\HH^m;
$

3) a curve $u(t)\in C^1(0,T;H)$ is a solution of \eqref{2.1} if and
only if $v(\tau)=\Psi(u(t))$ satisfies the equations
\begin{equation}\label{2.2}
\dot v_j=-{\rm sign}(j)\,v_{-j}\om_{|j|}(I_1,I_2,\dots)\,,
\quad j\in\Z_0\,,
\end{equation}
where $\om_l=\frac{\p \tilde K}{\p I_{l}}$ for $l=1,2,\dots$.
\end{theorem}

\begin{corollary}
If $u(t)$ is a solution of \eqref{2.1} and $\Psi(u)=v=\sum v_sf_s$,
then
\begin{equation}\label{2.3}
I_k(t)=\frac12\,(v_k^2+v_{-k}^2)(t)=\const\quad\forall\,k=1,2,\dots.
\end{equation}
\end{corollary}

If $v\in h^r$, then the vector $I=(I_1,I_2,\dots)$ belongs to the
space
$$
h_I^r=\{I:\, |I|_{h_I^r}=2\sum j^{1+2r}|I_j|<\infty\}.
$$
In fact, $I\in h^r_{I+}$, where
$$
 h^r_{I+}=
\{I\in h_I^r:\,
I_j\ge0\;\;\forall\,j\}\,.
$$

\noindent{\bf Amplification.} The function $\tilde K$ in
Theorem~\ref{t2.1} is analytic in $h^0_{I+}$. That is, it
analytically extends to the vicinity on this set in the space
$h^0_I$. \smallskip

The quantities $I_1,I_2,\dots$ are called {\it the actions\/}. Each
vector $\vv_j$ can be characterised by the action $I_j$ and {\it the
angle}
$$
\varphi_j=\arctan{  \frac{v_{-j}}{v_j} }\,.
$$
We will write $v=(I,\varphi)$, where
$\varphi=(\varphi_1,\varphi_2,\dots)$. The vector
$\varphi=(\varphi_1,\varphi_2,\dots)$ belongs to the torus
$\T^\infty$. We provide the latter with the Tikhonov topology, so it
becomes a compact set.

The functions $u\to v_k(u),\,k\in\Z_0$, form a coordinate system on
$H$. They are called the {\it Birkhoff coordinates}, and the system
of equations \eqref{2.2} -- the {\it Birkhoff normal form} for the
KdV equation. The normal forms is a classical tool to study
finite--dimensional Hamiltonian systems and their perturbations
locally in the vicinity of an equilibrium (see \cite{MoS}, \S30).
For all important finite--dimensional systems the normal forms do
not exist globally. In contrast, Theorem~\ref{t2.1} shows that the
KdV equation is an infinite--dimensional Hamiltonian system which
admits a normal form globally in the whole space $H$. To take all
advantages of this normal form we will need some information about
asymptotic properties of the transformation $\Psi(u)$ when
$u\to\infty$:

\begin{theorem}\label{t2.2}
For $m=0,1,\dots$ there are polynomials $P_m$ and $Q_m$ such that
$$
|d^j\Psi(u)|_m
\le P_m(\|u\|_m)\,, \quad j=0,1,2,
$$
and
$$
\|d^j\Psi^{-1}(v)\|_m
\le Q_m(|v|_m)\,,\quad j=0,1,
$$
for all $u,v$ and all $m\ge0$. Here for $j\ge1$ $|d^j\Psi|_m$ is the
norm of the corresponding poly-linear map from $H^m$ to $h^m$, and
similar with $\|d^j\Psi^{-1}\|_m$.
\end{theorem}

\begin{proof}
The estimates for the norms $|\Psi(u)|_m$ and $\|\Psi^{-1}(v)\|_m$
follows from Theorem~2.1 in \cite{Kor06}.\footnote{Note that the
quantity, denoted there $\|J\|_{p-1/2}$, equals $|v|_{p+1}$ up to a
constant factor, and $Q_{2p}$ satisfies the estimates $Q_{2p}\le
R_{1p}(\|u\|_{p+1})$ and $\|u\|_{p+1}\le R_{2p}(Q_{2p})$,  where
$R_{1p}$ and $R_{2p}$ are some polynomials.}

We do not prove here the  estimate for $d^j\Psi(u)$ with $j=1,2$. We
are certain that modern spectral techniques (e.g., see \cite{Kor06,
DM06})
 allow to establish them, but we think that this paper is not
a proper place for a corresponding rather technical research.
\end{proof}

\noindent{\bf Remark.} We do not use that the coordinate system
$v=(\vv_1,\vv_2,\dots)$ is symplectic, but only that it puts the KdV
equation to the form \eqref{2.2}. Therefore we may replace $v$ by
another smooth coordinate system $v'=(\vv'_1,\vv'_2,\dots)$ such
that $I_j'=I_j$ for all $j$ and
$\vp'_j=\vp_j+\Phi_j(I_1,I_2,\dots)$. Non-symplectic coordinate
systems are easier to construct, and it is possible that a proof of
Theorem~\ref{t2.2} simplifies if we replace there $v$ by a suitable
system $v'$.\qed\smallskip

For a function $f$ on a Hilbert space $H$ we write  $f\in\Plip(H)$
if
\begin{equation}\label{5.lip}
|f(u_1)-f(u_2)|\le P(R)\|u_1-u_2\|\;\;\text {if}\;\;\|u_1\|,\,
\|u_2\|\le R\,,
\end{equation}
where $P$ is a continuous function (depending on $f$). Clearly the
set of functions $\Plip(H)$ is an algebra. Due to the Cauchy
inequality any analytic function on $H$ belongs to $\Plip(H)$ (see
{\it Agreements}). In particular,
\begin{equation}\label{omega}
\om_l\in\Plip(h_r^I)\quad\text{for}\quad l\in\N,\, r\ge0.
\end{equation}


\section{Equation \eqref{1} in the Birkhoff coordinates}\label{s3}

For $k=1,2,\dots$ we denote
\begin{equation*}
    \Psi_k:H^m\to\R^2\,,\;\;\Psi_k(u)=\vv_k\,,
\end{equation*}
where $\Psi(u)=v=(\vv_1,\vv_2,\dots)$. Let $u(t)=u^\nu(t)$ be a
solution of \eqref{1}, which either is a stationary solution, or
satisfies \eqref{1.14} with a $\nu$-independent non--random $u_0$.
Applying Ito's formula to the map $\Psi_k$ we get:
\begin{equation}\label{3.1}
\begin{split}
d\vv_k=\big(d\Psi_k(u)(\nu u_{xx}+V(u))+&\frac12\,\nu\sum_{j\in\Z_0}
b_j^2
d^2\Psi_k(u)[e_j,e_j]\big)dt\\
 +&\sqrt\nu\,d\Psi_k(u)\big(\sum_{j\in\Z_0}
b_je_j\,d\beta^j\big).
\end{split}
\end{equation}
Let us denote
$$
d\Psi_k(u)\big(\sum b_je_j\,d\beta^j\big)= B_k(u)\,d\beta=
\sum_jB_{kj}(u)\,d\beta^j\,,\quad B_{kj}\in \R^2\;\;\forall\,k,j.
$$
Then the diffusion term in \eqref{3.1} may be written as
$\sqrt\nu\,B_k(u)\,d\beta$.

Since $I_k=\frac12\,|\Psi_k|^2$ is an integral of motion (see
\eqref{2.3}),  then application of Ito's formula to the functional
$\frac12\,|\vv_k|^2=I_k$ and eq.~\eqref{3.1} results in
\begin{equation}\label{3.3}
\begin{split}
 dI_k=\nu   &\Big((d\Psi_k(u){u_{xx}},\vv_k) +\frac12\,\big(\sum_j
 b_j^2d^2\Psi_k(u)[e_j,e_j],\vv_k\big)\\
 +&\frac12\,\sum_j
 b_j^2|d\Psi_k(u)e_j|^2\Big)dt+\sqrt\nu\,(B_k(u)\,d\beta,\vv_k)
\end{split}
\end{equation}
(here and below $(\cdot,\cdot)$ indicates the scalar product in
$\R^2$). Note that in difference with \eqref{3.1}, eq.~\eqref{3.3}
`depends only on the slow time' in the sense that all terms in its
r.h.s. have a factor $\nu$ or $\sqrt\nu$.

Let us consider the infinite-dimensional Ito process with components
\eqref{3.3}, $k\ge1$. The corresponding diffusion is
$\sqrt\nu\sigma\, d\beta$,  where   $\sigma=(\sigma_{kj}(u),
k\in\N,j\in\Z_0)$ and
$$
\sigma_{kj}= (B_{kj}(u),\vv_k)=
 \,b_j(d\Psi_k(u)e_j,\Psi_k(u)).
$$
Consider the {\it diffusion  matrix} $a$,
\begin{equation}\label{3.4}
a(u)=\sigma(u)\sigma^t(u)\,,\quad
a_{k_1k_2}=\sum_{j\in\Z_0}\sigma_{k_1j}\sigma_{k_2j}\,.
\end{equation}

\begin{lemma}\label{l3.1}
For any $u\in H$ the sums in \eqref{3.4} converge. The matrix $a$ is
symmetric and defines a bounded linear operator in $l^2$. If
$a\xi=0$ for some $\xi\in l^2$, then $\xi_k\ne0$ only if $\vv_k=0$,
where $v=\Psi(u)$. In particular, if $\vv_k\ne0\ \forall\,k$, then
{\rm Ker}$\, a=\{0\}$. Moreover, if $|v_j|\ge\delta$ for $1\le j\le
m$, then for any $\xi\in\R^m\times\{0\}\subset\R^\infty$ we have
\begin{equation}\label{3.est}
    \langle a(u)\xi,\xi\rangle_{l_2}=|\sigma^t(u)\xi|^2_{l_2}
    \ge C|\xi|_{l_2}^2\,,
\end{equation}
where $C$ depends on $\delta, m, |v|_1$ and the sequence $\{b_j\}$.
\end{lemma}
\begin{proof}
Using \eqref{4} and Theorem~\ref{2.1} we get that $|\sigma_{kj}|\le
C|j|^{-1}\eta_k$, where $\eta\in l^2$. Therefore $\sigma$ defines a
bounded linear operator $H\to l^2$ and $\sigma^t$ defines a bounded
operator $l^2\to H$. So $a=\sigma\sigma^t$ is a bounded operator in
$l^2$ and its matrix is well defined. Let us take any vector $\xi$.
Then $(a\xi,\xi)_{l^2}=\langle\sigma^t\xi,\sigma^t\xi\rangle$ where
\begin{equation}\label{4x}
(\sigma^t\xi)_j= \sum_kb_j(d\Psi_k(u)e_j,\vv_k)\xi_k= b_j\langle
e_j, d\Psi(u)^*(\oplus\xi_k\vv_k)\rangle.
\end{equation}
Hence, $\xi\in\,$Ker$\,a$ if and only if
 $d\Psi(u)^*(\oplus\xi_k\vv_k)=0$. Since $d\Psi(u)$ is an
isomorphism, then in this case $\xi_k\vv_k=0$ for each $k$, and the
assertion follows.

To prove \eqref{3.est} we abbreviate $\oplus\xi_k\vv_k=\xi_v$ and
denote $d\Psi(u)^*\xi_v=\eta$. Then
$\sigma^t(u)\xi=\,$diag$\,\{b_j\}\eta$ (see \eqref{4x}).
 Due to the first assertion
of Theorem~\ref{t2.2},
$$
\|\eta\|_1^2\le C_1(|v|_1)|\xi_v|_1^2\le C_1(|v|_1)\,|\xi|_{l_2}^2\,
C_2(|v|_0)m^3\,.
$$
So
$$
\sum_{k=N+1}^\infty \eta_k^2\le N^{-2}C_1C_2|\xi|^2_{l_2}m^3\,,
$$
for any $N$. Since $(d\Psi(u)^*)^{-1}=(d\Psi(u)^{-1})^*$,  then the
second assertion of the theorem implies that
$$
 \sum_{k=1}^\infty \eta_k^2=\|\eta\|_0^2\ge C_0(|v|_0)\,|\xi_v|^2_0\ge
C_0(|v|_0)\,|\xi|_{l_2}^2\delta^2\,.
$$
Choosing $N=
\big[\big(2C_1C_2C_0^{-1}\delta^{-2}m^3\big)^{1/2}\big]+1$ we get
that $\  \sum_{k=1}^N \eta_k^2\ge\frac12\,C_0\delta^2|\xi|^2_{l_2}
\,.$ Accordingly, $$  |\sigma^t(u)\xi|_{l_2}^2\ge C'\sum_{k=1}^N
\eta_k^2\ge \frac12\,C'C_0\delta^2|\xi|^2_{l_2}\,, $$ where $C'$
depends on the sequence $\{b_j\}$ and $N$.
\end{proof}

We see that the infinite--dimensional Ito process
\eqref{3.3}${}_{k\in\N}$, defined for $I\in h^0_{I+}$, has
non-degenerate diffusion outside the boundary $\p h^0_{I+}=\{I:\,
I_j=0$ for some $j\ge0\}$.

By applying Ito's formula to the $k$-th angle
$\displaystyle\varphi_k=\arctan\Big(\frac{v_{-k}}{v_k}\Big)$
$(k\ge1)$ and using \eqref{2.2} we obtain
\begin{equation*}\begin{split}
d\varphi_k&=\Big[ {\om}_{k}(I) + \nu|\vv_k|^{-2}(d\Psi_k(u)
u_{xx},\vv_k^{\bot})\\
+&\nu|\vv_k|^{-2}\big(\sum\limits_{j=1}^\infty b_j^2
d^2\Psi_k[e_j,e_j],\vv_k^{\bot}\big) -\nu|\vv_k|^{-2}
\sum_{j\in\Z_0}\big((B_{kj},\vv_k)(B_{kj},\vv_k^{\bot})\big)
\Big]dt\\
+&\sqrt{\nu}|\vv_k|^{-2}(B_k(u),\vv_k^{\bot})d\beta\,,
\end{split}
\end{equation*}
where $\vv_k^{\bot}= \left(\begin{array}{c}
  -  v_{-k}\\
    v_{k}\\
  \end{array}\right)$.
 Denote for brevity the drift and diffusion coefficients in the
above equation by ${\om}_k(I)+\nu G_k(v)$ and $\sqrt{\nu}\,g_k^j(v)$
respectively. Denoting similarly  the drift  coefficients in
(\ref{3.3}) by $\nu F_k(v)$ we rewrite the equation for the pair
$(I_k,\varphi_k)$ $(k\ge1)$ as
\begin{equation}\label{de_ug}\begin{split}
dI_k(t)&=\nu F_k(v)dt+\sqrt{\nu}\,\sigma_k(v)\,d\beta_t\,,
\\
d\varphi_k(t)&=[{\om}_k(I)+\nu
G_k(v)]dt+\sqrt{\nu}\,g_k(v)\,d\beta_t\,.
\end{split}
\end{equation}
\smallskip

Introducing the fast time  $$\tau=\nu t$$ we rewrite  the system
(\ref{de_ug}) as
\begin{equation}\label{de_ug_red}\begin{split}
dI_k(\tau)&= F_k(v)d\tau+\sigma_k(v)\,d\beta_\tau\,,
\\
 d\varphi_k(\tau)&=\Big[\frac{1}{\nu}\,{\om}_k(I)+
G_k(v)\Big] d\tau +g_k(v)\,d\beta_\tau\,.
\end{split}
\end{equation}
Here $\beta=(\beta_j,j\in\Z_0)$, where $\beta_j(\tau)$ are new
standard independent Wiener processes.

In the lemma below $P_k$ and $P_{kN}$ are some polynomials.

\begin{lemma}\label{l_polyzav} For $k\in\N, \,j\in\Z_0$ we have:

i) the function $F_k$  is analytic in each space $h^r$,  $r\ge2$ (so
$F_k\in\Plip(h^r)$), and has a polynomial growth as
$|v|_k\to\infty$;

ii) the function $\sigma_{kj}(v)$ is analytic in $h^r$, $r\ge0$,
 and for any $N\ge1$ satisfies $|\sigma_{kj}(v)|\le j^{-N}P_{kN}(|v|_r)\
 \forall\,v\in h^r$;

ii) for any $r\ge2,\delta>0$  and $N\ge1$ the functions
$G_k(v)\chi_{\{I_k>\delta\}}$ and $g_{kj}(v)\chi_{\{I_k>\delta\}}$
 are bounded, respectively, by $\delta^{-1}P_k(|v|_r)$ and
  $\delta^{-1}j^{-N}P_{kN}(|v|_r)$.
\end{lemma}
\begin{proof}
The assertions  concerning the functions $F_k$ and $G_k$ follow from
Theorem~\ref{t2.2} since the set of analytical functions with
polynomial growth at infinity is an algebra. To get the assertions
about $\sigma_k$ and $g_k$ we also use \eqref4.
\end{proof}

\section{More estimates}\label{s5}

In this Section and in the following  Sections~\ref{s.av}-\ref{s6}
we consider solutions of equation \eqref{de_ug}, written in the form
\eqref{de_ug_red}, which either are stationary in time, or satisfy
the $\nu$-independent initial condition \eqref{1.14}, where for
simplicity $u_0$ is smooth and non-random,
$$
u_0\in H^\infty=\bigcap_mH^m.
$$
 First we derive for these solutions additional
estimates, uniform in $\nu$.

\begin{lemma}\label{invmeabou}
For any $\nu>0,T>0 $ and $m,N\in\N$ the process $I(\tau)$ satisfies
the estimate
\begin{equation}\label{4.x}
\E\sup_{0\le \tau\le T}|I(\tau)|^N_{h^m_I}=\E\sup_{0\le \tau\le
T}|v(\tau)|^{2N}_{m}\le C(N,m,T)\,.
\end{equation}
\end{lemma}
\begin{proof}
For the sake of definiteness we consider a stationary solution
$v(\tau)=\{v_k^\nu(\tau)\}$. Cauchy problem
\eqref{de_ug_red},~\eqref{1.14} can be considered in the same way.
Applying Ito's formula to the expression $k^mI_k^N$ gives
$$
d(k^mI_k^N)=k^m\Big((NI_k^{N-1}F_k(v)+
\frac{1}{2}N(N-1)I_k^{N-2}\sum\limits_{j=1}^\infty(B_{kj}(v),\vv_k)^2)d\tau+
$$
$$
+NI_k^{N-1}\sigma_k(v)d\beta_\tau   \Big).
$$
Therefore,
$$
\E\sup_{0\le \tau\le T}k^mI_k^N(\tau)\le \E k^mI_k^N(0)+
$$
$$
+k^m\E\sup_{0\le \tau\le
T}\Big|\int\limits_0^{\tau}\Big(NI_k^{N-1}(s)F_k(v)+
\frac{1}{2}N(N-1)I_k^{N-2}(s)\sum\limits_{j=1}^\infty\sigma_{kj}^2\Big)\,ds
\Big|+
$$
$$
k^m\E\sup_{0\le \tau\le
T}\Big|\int\limits_0^{\tau}NI_k^{N-1}(s)\sigma_k(v)d\beta_s\Big|\le
 C(m,N,T)\,.
$$
Doob's inequality, Lemma~\ref{l_polyzav} and Theorem~\ref{tI.1} have
been used here. This relation yields the desired estimate. Indeed,
by the H\"older inequality we get
\begin{equation*}\begin{split}
&\E(\sup_{0\le \tau\le T}|I(\tau)|_m^{2N})= 
2^N \E\sup_{0\le \tau \le
T}\!\Big(\sum\limits_{j=1}^\infty\frac{1}{j^2}
j^{2m+3}I_j(\tau)\Big)^N\\
&\le2^N \E\sup_{0\le \tau\le T}\left\{\Big(\sum\limits_{j=1}^\infty
j^{N(2m+3)}I_j^{N}(\tau)\Big)^{N\frac{1}{N}}\Big(\sum\limits_{j=1}^\infty
j^{-\frac{2N}{N-1}}\Big)^{N\frac{N-1}{N}}\right\}\\
&\le C_N\E\sup_{0\le \tau\le T}\Big(\sum\limits_{j=1}^\infty
j^{N(2m+3)}I_j^{N}(\tau)\Big)\le 
C_1(m,N,T).
\end{split}
\end{equation*}
\end{proof}

In the further analysis we  systematically use the fact that the
functionals $F_k(I,\varphi)$ depend weakly on the tails of
vectors $\varphi=(\varphi_1,\varphi_2,\dots$). 
Now we state the corresponding auxiliary results.

Let $f\in\Plip(h^{n_1})$ and $v\in h^n$, $n>n_1$. Denoting by
$\Pi_M$, $M\ge1$, the projection
$$
\Pi_M:h^0\to h^0,\quad \sum v_jf_j\mapsto\sum_{|j|\le M}v_jf_j\,,
$$
we have $\  |v-\Pi_Mv|_{n_1}\le M^{-(n-n_1)}|u|_n$. Accordingly,
\begin{equation}\label{5.new}
|f(v)-f(\Pi_M(v))|\le P(|v|_n)M^{-(n-n_1)}.
\end{equation}
Similar inequalities hold for functions on $h^n_I$, and
\eqref{omega} with $r=0$ implies that
\begin{equation}\label{5.om}
|\om_k(I)-\om_k(\Pi_MI)|\le P_k(|I|_n)M^{-n}.
\end{equation}

The torus $\T^M$ acts on the space $\Pi_Mh^0$ by linear
transformations $\Phi_{\theta_M}$, $\theta_M\in\T^M$, where
$\Phi_{\theta_M}$ sends a point $v_M=(I_M,\varphi_M)$ to
$(I_M,\varphi_M+\theta_M)$. Similar, the torus $\T^\infty$ acts on
$h^0$ by linear transformations
$\Phi_\theta:(I,\varphi)\mapsto(I,\varphi+\theta)$. The
transformation $\Phi_\theta$ continuously depends on
$\theta\in\T^\infty$, in the strong operator topology.

For a function $f\in\Plip(h^{n_1})$ and any $N$ we define the
average of $f$ in the first $N$ angles as the function
$$
\langle
f\rangle_N(v)=\int_{\T^N}f\big(\Phi_{\theta_N}\oplus\text{id}\,)
(v)\big)\,d\theta_N
$$
(here $id$ stands for the identity transformation in the space
$h^0\ominus\Pi_Nh^0$), and define the average in all angles as
$$
\langle f\rangle(v)=\int_{\T^\infty}f(\Phi_\theta v)\,d\theta\,,
$$
where $d\theta$ is the Haar measure on $\T^\infty$. The estimate
\eqref{5.new} readily implies that
\begin{equation}\label{5.est}
|\langle f\rangle_N(v)-\langle f\rangle(v)| \le
P(R)N^{-(n-n_1)}\;\;\text{if}\;\,|v|_n\le R\,.
\end{equation}

Let $v=(I,\varphi)$. Then $\langle f\rangle_N$ is a function,
independent of $\varphi_1,\dots,\varphi_N$, and $\langle f\rangle$
is independent of $\varphi$. I.e., $\langle f\rangle$ can be written
as a function $\langle f\rangle(I)$.

\begin{lemma}\label{l_aver} Let $f\in \Plip(h^{n_1})$. Then

i) The functions $\langle f\rangle_N(v)$ and $\langle f\rangle(v)$
satisfy \eqref{5.lip} with the same polynomial as $f$ and take the
same value at the origin.

ii) They are smooth (analytic) if $f$ is. Moreover, if $f$ is
smooth, then $\langle f\rangle(I)$ is a smooth functions of the
vector $(I_1,\dots,I_M)$ for any $M$. If $f(v)$ is analytic in the
space $h^{n_1}$, then $\langle f\rangle(I)$ is analytic in the space
$h_I^{n_1}$.

\end{lemma}
\begin{proof} i) Is obvious.

ii) The first assertion is obvious. To prove the last two consider
the function $g(r_1,r_2,\dots)=\langle f\rangle(\vv_1,\vv_2,\dots)$,
$\vv_j=\left(
         \begin{array}{c}
           r_j \\
           0 \\
         \end{array}
       \right)$.
Then $g(r)=\langle f\rangle(I)$, where $I_l=\frac12r_l^2$ for each
$l$. The function $g$ is smooth and even in each $r_j$, $j\ge1$. Any
function of finitely many arguments with this property is known to
be a smooth function of the squared arguments, so the second
assertion holds.

Now let $f(v)$ be analytic. Denote by ${\frak h}^{n_1}$ the space of all
sequences $r=(r_1,r_2,\dots)$ such that the corresponding vector $v$ belongs to
$h^{n_1}$, and provide it with the natural norm. If $f(v)$ is analytic, then
$\langle f\rangle(v)$ also is analytic and
$g(r)$ extends analytically to an even function in
a complex neighbourhood $\cal O$ of ${\frak
h}^{n_1}$ in ${\frak h}^{n_1}\otimes\C$. This neighbourhood may be chosen to be
invariant with respect to all involutions
$$
(r_1,r_2,\dots,r_j,\dots)\mapsto (r_1,r_2,\dots,-r_j,\dots), \quad
j=1,2,\dots.
$$
The image ${\cal O}_I$  of $\cal O$ under the map
$$
(r_1,r_2,\dots)\mapsto (\tfrac12 r_1^2,\tfrac12r_2^2,\dots)
$$
is a neighbourhood of $h_I^{n_1}$ in the complex space
$h_I^{n_1}\otimes\C$.  The function
$$
g(\pm\sqrt{2I_1},\pm\sqrt{2I_2},\dots)=:g(\sqrt{I})
$$
is a well defined locally  bounded function on   ${\cal
O}_I$.\footnote{i.e., it is bounded uniformly on bounded subsets of
${\cal O}_I$.} For any $N$ its restriction to  ${\cal O}_I^N= {\cal
O}_I \cap\Pi_N(h_I^{n_1}\otimes C)$ is a single-valued algebraic
function on a domain in $\C^N$; so $g(\sqrt{I})$ is analytic on
${\cal O}_I^N$ for each $N$. Hence, $g(\sqrt{I})$ is analytic on
${\cal O}_I$ (see Lemma~A.4 in \cite{KaP}). Since
$g(\sqrt{I})=\langle f\rangle(I)$, then the result follows.
\end{proof}

Let $(I^\nu(\tau),\varphi^\nu(\tau))$ be a solution of
\eqref{de_ug_red}.
 In the lemma below we show that the processes
$I_k^\nu(\tau),\,k\ge1$, do not asymptotically approach zero as
$\nu\to 0$ (concerning the notation $\vk(\delta^{-1};M,T)$, used
there, see {\it Notations\/}):
\begin{lemma}\label{notatzero}
For any $M\in\mathbb N$ and $T>0$ we have
\begin{equation}\label{4.y}
{\bf P}\{\mathop{\rm min}\limits_{k\le
M}I_k^{\nu}(\tau)<\delta\}\le\vk(\delta^{-1};M,T),
\end{equation}
uniformly in  $\nu>0$ and $0\le\tau\le T$.
\end{lemma}

Here the difficulty is that the scalar process
$I_k(\tau)=\frac12|\vv_k(\tau)|^2$ satisfies equation
\eqref{de_ug_red}, where the diffusion $\sigma_k$ degenerates when
$I_k$ vanishes. The equation for the vector-process $\vv_k(\tau)$
(see \eqref{scal3.1} below) has a non-degenerate diffusion, but its
drift has a component of order $\nu^{-1}$. To prove the lemma's
assertion we construct a new process $\hat\vv_k(\tau)$ such that
$|\vv_k(\tau)|=|\hat\vv_k(\tau)|$ and $\hat\vv_k$ satisfies an Ito
equation with a nondegenerate diffusion and coefficients, bounded
uniformly in $\nu$. Then $I_k=\frac12|\hat\vv_k(\tau)|^2$ meets
estimate \eqref{4.y} by a Krylov's theorem. The problem to perform
this scheme is that the process $\hat \vv_k$ is constructed as a
solution of an additional diffusion equation which is ill defined
when $v_k$ vanishes. We cannot show that the event
$$
\{v_k(\tau)=0 \text { for some } 0\le \tau\le T\}
$$
has zero probability and resolve this new difficulty by means of
some additional (rather involved) construction.

For a complete proof see  Section \ref{proof}.


\section{Averaging along Kronecker flows.}\label{s.av}

The flow
$$
S^t:\T^\infty\to\T^\infty\,,\quad \varphi\mapsto \varphi+tW,\;\;
t\in\R\,,
$$
where $W\in \R^\infty$, is called a {\it Kronecker flow}. In this
section we study averages of functions $f(v)=f(I,\varphi)$ along
such flows. That is, we study the quantities
$$
\frac{1}{T}\int\limits_0^T f(I,\varphi+{\om}^m t)dt,
\quad T>0\,.
$$

\begin{lemma}\label{proobraz}
Let $f\in\Plip(h^{n_1}), v=(I,\varphi)\in h^n, n>n_1\ge0$, and $f$ is
analytic in the space $h^{n_1}$. Then for  each ${R'}>0$,
$m\in\mathbb N$ and $\delta>0$ there is a Borel  set
${\Omega}_{R'}^m(\delta)\subset\{x\in\mathbb R^m\,:\, |x|\le {R'}\}$
such that $|{\Omega}_{R'}^m(\delta)|<\delta$, and for any
${\om}^m\not\in {\Omega}_{R'}^m(\delta)$, $|\om^m|\le {R'}$ the
estimate
$$
\left|\frac{1}{T}\int\limits_0^T f(I,\varphi+{\om}^m t)dt- \langle
f\rangle(v)\right|
 \le\frac{1}{T\delta}c_0(m,{R'},|v|_n,  f)+m^{-(n-n_1)}P(|v|_n)\,,
$$
holds uniformly in $\varphi\in\mathbb T^\infty$. Here $P$ is the
continuous function from \eqref{5.lip} and $W^m$ is identified with
the vector $(W^m,0,\dots)\in\R^\infty$.
\end{lemma}
\begin{proof}
Let us first assume that   $f(v)=f(\Pi_mv)$ (i.e.,
 $v$ depends only on finitely-many
variables). Then $f=f(I^m,\varphi^m)$ is analytic in $\varphi^m$ and
 the radius of analyticity is independent of $I$, satisfying
 $|I|_{h^0_I}\le R'$. Now the estimate with $P:=0$ is a classical result
(e.g., see in \cite{MoS}). In general case we write $f$ as
$f\circ\Pi_m+ (f-f\circ\Pi_m)$ and use \eqref{5.est}.
\end{proof}
\medskip

We will apply this lemma with $W^m=W^m(I)$, where $I=I(\tau)$ is the
$I$-component of a solution of \eqref{de_ug}. To do this we have to
estimate probabilities of the events $\{W^m(I(\tau))\in
\Omega^m_{R'}(\delta)\}$. To state the corresponding result we
introduce more notations. For any events $Q$ and $\cal O$ we denote
$$
\IP_{Q}({\cal O})=\IP\big(\complement Q \cap\cal O\big)\,,
$$
and
$$
\E_{Q}(f)=\E\big((1-\chi_{Q})  f\ \big)\,.
$$
Abusing language, we call $\IP_Q$ a probability.  We fix any
$$
p\ge1,
$$
denote
$$
\BB_R=\{I:\, |I|_{h_I^p}\le R\},
$$
and for $R\ge1$ consider the event
$$
\Omega_R=\{\sup_{0\le\tau\le T}|v^\nu(\tau)|_p\ge R\}\,,
$$
where $v^\nu(\tau)$ is a solution. Noting that $|W^m(I)|\le
R'=R'(R,m)$ outside the event $\Omega_R$, we denote
$$
\Omega(\delta)=\Omega^m_{R'}(\delta),\quad R'=R'(R)\,,\
0<\delta<1\,.
$$
Finally, for $M\ge m$ and $0<\gamma<1$ we define
$$
Q_\gamma=\{I\in h^0_{I+}:\, \min_{1\le j\le M}I_j<\gamma\}\,.
$$

\begin{lemma}\label{l.new}
There exists $M=M(R,m)\ge m$ such that
\begin{equation}\label{y0}
    \int_0^T\PP_{\Omega_R}\Big(\{W^m(I(s))\in\Omega(\delta)\}
    \setminus\{I(s)\in Q_\gamma\}\Big)\,ds
    \le\vk(\delta^{-1};R,m,\gamma,T)\,,
\end{equation}
uniformly in $\nu>0$. \footnote{We recall that $\vk(t;R,m,\gamma,T)$
stands for a function of $t$ which goes to zero when $t\to\infty$,
and depends on the parameters $R,m,\gamma$ and $T$.}
\end{lemma}

\begin{proof}
Consider the function $D(I)=\det\big(\p W^m_j/\p I_r: 1\le j,r\le
m\big)$. It is analytic in $h_I^0$ (see Amplification to
Theorem~\ref{t2.1}), and $D\not\equiv0$ since $D(0)=C^m$, $C\ne0$
(see \cite{K2}, Lemma~3.3,   and \cite{KaP}).  For a finite
non-decreasing sequence of
natural numbers $\alpha=(\alpha^1\le\dots\le\alpha^N)$ we denote
$$
|\alpha|=\alpha^N\,,\quad [\alpha]=N
$$
and define the derivative $\p^\alpha D(I)/\p I^\alpha$ in the
natural way.\smallskip

{\bf Step 1:} Study of the sets
$\{I\in\BB_R:\,|D(I)|<\e\},\;0<\e\ll1$.

By the analyticity  any point $I'\in\BB_R$ has a neighbourhood
${\cal O}\subset h^0_I$ such that
$$
\Big|\frac{\p^\alpha D(I)}{\p I^\alpha}\Big|\ge
c\quad\forall\,I\in\cal O\,,
$$
where the sequence $\alpha=(\alpha^1\le\dots\le\alpha^N)$ and $c>0$
depend only on the neighbourhood. Since $\BB_R$ is a compact subset
of $h^0_I$, we can cover it by a finite system of neighbourhoods
${\cal O}_j$, $j=1,\dots,L$, as above, where $L=L(R,m)$. Then
\begin{equation}\label{y1}
    \{I\in\BB_R:\,\Big|\frac{\p^{\alpha_j}D(I)}{\p
    I^\alpha_j}\Big|\ll1\,, \quad j=1,\dots,L\}=\emptyset.
\end{equation}

Let us denote
$$
M=\max_{1\le j\le L}|\alpha_j|\,,\;\;N=\max_{1\le j\le L}[\alpha_j]
$$
 and consider the sequence
$$
\e=\e_0<\e_1<\dots<\e_N<1\,,\;\;\e_j=\e^{2^{-j}-2^{-N}+2^{-j-N}}\,,
$$
where $0<\e<1$. Note that
$$
\e_j\,\e^{-2}_{j+1}=\e^{(2^{-N})}\quad\text{for}\quad 0\le j<N\,.
$$

For $m\le[\alpha_j]$ we set
$$
\fA^m=\{I\in\BB_R:\,\Big|\frac{\p}{\p\alpha_j^1}
\dots\frac{\p}{\p\alpha_j^m}D(I)\Big|<\e_m\}\,.
$$
In particular, $\ \fA^0={\frak A}^0=\{I\in\BB_R:\, |D(I)|\le\e\}$
for each $j$.

For $0<\e\ll1$ relation \eqref{y1} implies that
$$
{\frak A}^0=\bigcup_{j=1}^L\Big(\big({\frak A}^0\setminus
\fA^1\big)\cup
\big(\fA^1\setminus\fA^2\big)\cup\dots\cup\big(\fA^{[\alpha_j]-1}\setminus
\fA^{[\alpha_j]}\big)\Big)\,.
$$\smallskip

{\bf Step 2:} An estimate for the integral
$\int_0^T\PP_{\Omega_R}\{|D(I(s))|<\e\}\,ds$.

Due to the last displayed formula, the integral to be  estimated is
bounded by a finite sum of the terms
\begin{equation}\label{y2}
    \int_0^T \PP_{\Omega_R}\{I(s)\in \fA^r\setminus\fA^{r+1}\}
    \,ds\,,\quad r< [\alpha_j]\,.
\end{equation}
To estimate \eqref{y2}, we abbreviate $\frac{\p}{\p\alpha_j^1}
\dots\frac{\p}{\p\alpha_j^r}D(I)=f(I)$. Then
\begin{equation}\label{ho}
\fA^r\setminus\fA^{r+1}=\{I\in\BB_R:\, |f(I)|<\e_r\;\;\text{and}\;\;
\Big|\frac{\partial}{\partial (\alpha_j^{r+1})}f(I)\Big|\ge
\e_{r+1}\}\,.
\end{equation}
Consider the Ito process
$z(\tau)=f(I(\tau))$. We define the Markov moment %
$\tau'=\min\{\tau\ge0:\,|I(\tau)|_{h^p_I}\ge R^2\}\wedge T$, and
re-define $z(\tau)$ for $\tau\ge\tau'$ as a continuous process,
satisfying
$$
dz(\tau)=d\beta^1_\tau\quad\text{for}\quad \tau\ge\tau'\,.
$$
Since $\tau'>T$ outside $\Omega_R$, then outside $\Omega_R$ we have
  $z(\tau)=f(I(\tau))$ for
$0\le\tau\le T$. For $z(\tau)$ we have
$$
dz(\tau)=c(\tau)\,d\tau+\sum b_j(\tau)\,d\beta^j_\tau\,,
$$
where $|c|\le C(R,m)$, $b_j=\delta_{j,1}$ for $\tau\ge\tau'$ and
$b_j=\sum\frac{\p f}{\p I_k}\sigma_{kj}$ for $\tau\le\tau'$.
Denoting $a=\sum b_j^2$, we have
$a=\sum(\sigma\sigma^t)_{jk}\nabla_jf\nabla_kf$. So $|a(\tau)|\le
C(R,m)$. From other hand, \eqref{3.est} in Lemma~\ref{l3.1} implies
that
\begin{equation}\label{y3}
    |a(\tau)|\ge
    C(R,m,\gamma)\,\sum_{j=1}^{M}(\nabla_jf)^2\quad\text{if} \quad
    I(\tau)\notin Q_\gamma\,.
\end{equation}

Applying Theorem~2.3.3 from \cite{Kr77} to the process $z(\tau)$, we
get
$$
\E\int_0^T(\chi_{\{|z(\tau)|\le\e_r\}}|a(\tau)|\,d\tau\le
C(R,m,T)\e_r\,.
$$
By \eqref{ho} and \eqref{y3} the integrand is
$\ge\e_{r+1}^2C(R,m,\gamma)$ if
$I(t)\in(\fA^r\setminus\fA^{r+1})\setminus Q_\gamma$. Hence,
\begin{equation*}\label{y4}
\begin{split}
\int_0^T\PP_{\Omega_R}\{I(s)\in(\fA^r\setminus\fA^{r+1})\setminus
Q_\gamma\}\,ds&\le \e_r\e_{r+1}^{-2}C(R,m,\gamma,T)\\
&= \e^{(2^{-N})}C(R,m,\gamma,T)\,.
\end{split}
\end{equation*}
We have seen that
\begin{equation}\label{y5}
\int_0^T\PP_{\Omega_R}\big(\{|D(I(s))|<\e\}\setminus \{I(s)\in
Q_\gamma\}\big)\,ds\le\e^{(2^{-N})}C_1(R,m,\gamma,T)\,.
\end{equation}\smallskip

{\bf Step 3:} Proof of \eqref{y0}.

We have an inclusion of events
\begin{equation*}\begin{split}
\{W^m(s)\in\Omega(\delta)\}\setminus\{I(s)\in Q_\gamma\}
&\subset\Big[
\Big(\{W^m(s)\in\Omega(\delta)\}\setminus\big(\{I(s)\in Q_\gamma\}\\
\cup\{D(I(s))<\e\}\big)\Big)
&\cup\big(\{|D(I(s))<\e\}\setminus\{I(s)\in Q_\gamma\}\big)\Big]\,.
\end{split}
\end{equation*}
Probability of the second event in the r.h.s. is already estimated.
To estimate probability of the first event we apply
the Krylov estimate to the process $W^m(s)$. Re-defining it after
the moment $\tau'$ (see Step~2) and arguing as when deriving
\eqref{y5} we get that
\begin{equation}\label{y6}\begin{split}
    \int_0^T \PP_{\Omega_R}\Big(\{W^m(s)\in\Omega(\delta)\}&
    \setminus\big(\{I(s)\in Q_\gamma\}
\cup\{D(I(s))<\e\}\big)\Big)\,ds\\
&\le|\Omega(\delta)|^{1/m}C(R,m,\gamma,\e,T)\,.
\end{split}
\end{equation}
Finally, choosing first $\e$ so small that the r.h.s. of \eqref{y5} is
$\le\tilde\e$ and next choosing $\delta$ so small that the r.h.s. of
\eqref{y6} is $\le\tilde\e$, we see that the l.h.s. of \eqref{y0} is
$\le2\tilde\e$ for any $\tilde\e>0$, if $\delta$ is sufficiently
small.
\end{proof}

\section{The limiting dynamics.}\label{s6}

Let us fix any $T>0$, an integer $p\ge3$ and abbreviate
$$
h^p=h\,,\;\;h_I^p=h_I\,,\;\;h^p_{I+}=h_{I+}\,\;
|I|_{h^p_I}=|I|\,,\,\;\; |v|_p=|v|\,.
$$

 Due to Lemma~\ref{invmeabou} and the equation,
satisfied by $I^\nu(\tau)$, the laws $ {\cal L}\{I^\nu(\cdot)\}\  $
form a tight family of Borel measures on the space
$C([0,T];h_{I+})$. Let us denote by $\QQ^0 $ any its weak limiting
point:
\begin{equation}\label{6.0}
\QQ^0=\lim\limits_{\nu_j\to 0}{\cal L}\{I^{\nu_j}(\cdot)\}\,.
\end{equation} Our aim is to show that $\QQ^0$ is a solution to the
martingale problem in the space $h_I$ with the drift operator
$\langle F\rangle(I)=(\langle F_1\rangle(I),\langle
F_2\rangle(I),\dots)$ and the covariance $\langle A\rangle(I)
=\{\langle A_{kl}\rangle(I)\}$, where
$$
\langle A_{kl}\rangle (I)=
\big\langle\big(\sigma(v)\sigma^t(v)\big)_{kl} \big\rangle=
\big\langle\sum\limits_jb_j^2\big(d\Psi_k(u)e_j,\vv_k\big)
\big(d\Psi_l(u)e_j,\vv_l\big) \big\rangle.
$$
By Lemmas~\ref{l_polyzav} and \ref{l_aver} the averages $\langle
F_j\rangle$ and $\langle A_{kl}\rangle$ are analytic functions on
$h_I$. The covariance $\langle A\rangle$ is non-degenerate outside
the boundary of the domain $h^p_{I+}$ in the following sense: let
$\xi\in\R^M\subset\R^\infty$ and $I\in h_{I+}, |I|\le R$.
 Then
\begin{equation}\label{6.cov}
    \sum_{k,l\le M}\langle A_{kl}\rangle (I)\xi_k\xi_l \ge
    C|\xi|^2_{l_2}\;\;\;\text{if}\;\;\;
    |I_j|\ge\gamma>0\;\;\text{for}\;\;j\le M\,,
\end{equation}
where $C>0$ depends on $M, R$ and $\gamma$. Indeed, the estimate
follows from \eqref{3.est} with $v=(I,\varphi)$ by averaging in
$\varphi$.

Our study of the limit $\QQ^0$ uses the scheme, suggested by
R.~Khasminskii in \cite{Khas68} and is  heavily based on the
estimates for solutions $v^\nu(\tau)$, obtained above.

First we show that for any $k$ the difference
\begin{equation}\label{8.4}
I_k(\tau)-\int\limits_0^{\tau} \langle F_k\rangle(I(s))ds
\end{equation}
is a martingale with respect to $\QQ^0$ and the natural filtration
of $\sigma$-algebras.  A crucial step of the proof is to establish that
\begin{equation}\label{driftlim}
{\frak A}^\nu:= {\bf E}\,\max\limits_{0\le \tau\le
T}\bigg|\int\limits_0^{\tau}\big(F_k(I^\nu(s),
\varphi^\nu(s))-\langle F_k\rangle(I^\nu(s))\big)\,ds\bigg|\,\to\,0
\end{equation}
 as $\nu\to0$.  Proof of \eqref{driftlim} occupies most of this section.

 Let us fix an integer
 $$
  m\ge1,
 $$
  denote the first $m$ components of vectors $I^\nu$  and $\varphi^\nu$
 by $I^{\nu,m}$  and $\varphi^{\nu,m}$, and
rewrite the first $2m$ equations of  the system (\ref{de_ug_red}) as
follows
\begin{equation}\label{findims}
\begin{array}{rcl}\displaystyle
dI^{\nu,m}&\!\!=\!\!&
F^m(I^\nu,\varphi^\nu)d\tau+\sigma^m(I^\nu,\varphi^\nu)
d\beta_\tau\,,\\[2mm]
\displaystyle
d\varphi^{\nu,m}&\!\!=\!\!&\displaystyle\Big(\frac{1}{\nu}
{\om}^m(I^{\nu})+
G^m(I^\nu,\varphi^\nu)\Big)\,d\tau
 +g^m(I^\nu,\varphi^\nu)d\beta_\tau\,.
\end{array}
\end{equation}
Here and afterwards  we identify the vectors
$(I^\nu_1,\dots,I^\nu_m,0,0,\dots)$ with $I^{\nu,m}$, and the
vectors $(\varphi^\nu_1,\dots,\varphi^\nu_m,0,0,\dots)$ with
$\varphi^{\nu,m}$.

Denote $\  \langle F_k\rangle_m(I^m)= \langle
F_k\rangle_m(I,\varphi)_{I=(I_m,0),\varphi=0}$.
 By  Lemma~\ref{l_polyzav} 
there is a constant $C_k(R)$ such that for any $v=(I,\varphi)$,
$|v|\le R$,  we have
\begin{equation}\label{ine_1}
|F_k(I,\varphi)- F_k(I^{m},\varphi^{m})|\le C_k(R)m^{-1},
\end{equation}
\begin{equation}\label{ine_2}
|\langle F_k\rangle_m(I^m)-\langle F_k\rangle(I) |\le
C_k(R)m^{-1}\,.
\end{equation}
%

Define the event $\Omega_R$ as in Section~\ref{s.av}. Due to
Lemma~\ref{invmeabou}
$$
\PP(\Omega_R)\le\vk_\infty(R)
$$
(here and in similar situations below the function $\vk$ is
$\nu$-independent). Since by Lemma~\ref{l_polyzav} the function
$F_k$ has a polynomial growth in $v$, then this estimate implies
that
$$
|\,\E\max_{0\le\tau\le T}\int_0^\tau F_k(v^\nu(s))\,ds-
\E_{\Omega_R}\max_{0\le\tau\le T}\int_0^\tau F_k(v^\nu(s))\,ds|\le
\varkappa_\infty(R).
$$
 The functions
$F_k(I^{\nu,m},\varphi^{\nu,m})$, $\langle F_k\rangle_m(I^{\nu,m})$
and $\langle F_k\rangle(I^{\nu,m})$ satisfy similar relations. So we
have
$$
{\frak A}^\nu\le\ \varkappa_\infty(R)+
\E_{\Omega_R}\,\max\limits_{0\le \tau\le
T}\left|\int\limits_0^{\tau}\{F_k(I^\nu(s),\varphi^\nu(s))ds-
F_k(I^{\nu,m}(s),\varphi^{\nu,m}(s))\}ds\right|
$$
$$
+\E_{\Omega_R}\,\max\limits_{0\le \tau\le
T}\left|\int\limits_0^{\tau}\{F_k(I^{\nu,m}(s),\varphi^{\nu,m}(s))-
\langle F_k\rangle_m(I^{\nu,m}(s))\}ds\right|\,
$$
$$
+\E_{\Omega_R}\,\max\limits_{0\le \tau\le
T}\left|\int\limits_0^{\tau}\{\langle F_k\rangle_m(I^{\nu,m}(s))-
 \langle F_k\rangle(I^{\nu}(s))\}ds\right|\,\le
$$
\begin{equation*}
\begin{split}
\le \varkappa_\infty(R)+&  C_k(R)m^{-1} \\
&+\E_{\Omega_R}\,\max\limits_{0\le \tau\le
T}\left|\int\limits_0^{\tau}\{F_k(I^{\nu,m}(s),\varphi^{\nu,m}(s)) -
\langle F_k\rangle_m(I^{\nu,m}(s))\}ds\right|\,.
\end{split}
\end{equation*}
The last inequality here follows from (\ref{ine_1})-(\ref{ine_2}).
It remains to estimate the quantity
$$
\max\limits_{0\le \tau\le T}
\left|\int\limits_0^{\tau}\{F_k(I^{\nu,m}(s),\varphi^{\nu,m}(s))-
\langle F_k\rangle_m(I^{\nu,m}(s))\}ds\right| (1-\chi_{\Omega_R})\,.
$$
To do this  we consider a  partition of the interval $[0,T]$ to
subintervals of length $\nu L, L>1$ by the points
$$
\tau_j=\nu t_0+\nu jL,\quad 0\le j\le K+1 \,,
$$
where $\tau_{K+1}$ is the last point $\tau_j$ in $[0,T]$.  The
constant $L$ such that
\begin{equation}\label{apr1}
L\ge2,\quad L\le \frac12\, \nu^{-1}
\end{equation}
and the  (deterministic) initial point  $t_0\in[0,L)$
  will be chosen later. Note that
$$
\frac12\, T\le K\cdot \nu L\le T.
$$

 Denote
$$
\eta_l=\int\limits_{\tau_l}^{\tau_{l+1}}
\Big(F_k(I^{\nu,m}(s),\varphi^{\nu,m}(s))ds- \langle
F_k\rangle_m(I^{\nu,m}(s))\Big)ds\,,\quad 0\le l\le K\,.
$$
Since outside the event $\Omega_R$ we have
$$
\bigg |\int\limits_{\tau'}^{\tau''}
\Big(F_k(I^{\nu,m}(s),\varphi^{\nu,m}(s))- \langle
F_k\rangle_m(I^{\nu,m}(s))\Big)ds\bigg| \le \nu LC(R)
$$
for any $\tau'<\tau''$ such that $\tau''-\tau'\le\nu L$, then
\begin{equation}\label{6x}
\begin{split}
\E_{\Omega_R}\max\limits_{0\le \tau\le T}\bigg|&\int\limits_{0}^\tau
\big(F_k(I^{\nu,m}(s),\varphi^{\nu,m}(s))- \langle
F_k\rangle_m(I^{\nu,m}(s))\big)ds\bigg|\\
&\le\E_{\Omega_R}\, \sum\limits_{l=0}^{K}|\eta_l|+\nu LC(R)\,.
\end{split}\end{equation}
 To calculate the contribution from the integral over an $l$-th
subinterval, we pass there to the slow time $t=\nu^{-1}\tau$. Now
the system (\ref{findims}) reads as
\begin{equation}\label{slowfindims}
\begin{array}{rcl}\displaystyle
dI^{\nu,m}(t)&\!\!=\!\!& \nu
F^m(I^\nu,\varphi^\nu)dt+\sqrt{\nu}\sigma^m(I^\nu,\varphi^\nu)
d\beta_t\,,\\[2mm]
\displaystyle d\varphi^{\nu,m}(t)&\!\!=\!\!&\displaystyle\big(
{\om}^m(I^{\nu}) 
+\nu G^m(I^\nu,\varphi^\nu)\big)\,dt
 +\sqrt{\nu}g^m(I^\nu,\varphi^\nu)d\beta_t\,.
\end{array}
\end{equation}

Denoting $t_j=\tau_j/\nu=t_0+jL$ we  have:
$$
|\eta_l|\le \nu\bigg|\!\!\int\limits_{t_l}^{t_{l+1}}\!\!\!
\big\{F_k(I^{\nu,m}(x),\varphi^{\nu,m}(x))-\hskip6cm
$$
$$
\hskip1cm -F_k\Big(I^{\nu,m}(t_l),\varphi^{\nu,m}(t_l)+
{\om}^m(I^{\nu}(t_l))(x-t_l)\Big)\big\}dx\bigg|
$$
$$
+\nu\bigg|\int\limits_{t_l}^{t_{l+1}} \big\{F_k\Big(I^{\nu,m}
(t_l),\varphi^{\nu,m}(t_l)+{\om}^m(I^{\nu}(t_l))(x-t_l)\Big)
 -\langle F_k\rangle_m(I^{\nu,m}(t_l))\big\}dx\bigg|
$$
$$
+\nu\bigg|\int\limits_{t_l}^{t_{l+1}} \big\{ \langle
F_k\rangle_m(I^{\nu,m}(t_l))- \langle
F_k\rangle_m(I^{\nu,m}(x))\big\}dx\bigg|\,=\Upsilon_l^1+\Upsilon_l^2+
\Upsilon_l^3\,.
$$

To estimate the integrals $\Upsilon^1_l-\Upsilon^3_l$ we first
optimise the choice of $t_0$. Defining the event $\Omega(\delta)$,
the number $M(R,m)$ and the set $Q_\gamma$ as in Section~\ref{s.av},
we have

\begin{lemma}\label{lmeas}
The non-random number $t_0\in[0,\nu L)$ (depending on $\nu$ and $\delta$)
can be chosen in such a way that
\begin{equation}
  \label{meas}
  \frac1{K}\,\sum_{l=0}^K\PP{\cal E}_l\le\vk_\infty(R)+\vk(\gamma^{-1};R,m)
  +\vk(\delta^{-1};\gamma,R,m)\
\end{equation}
for all $0<\delta,\gamma<1$, where
$$
\cE_l=\Omega_R\cup \{I(\tau_l)\in
Q_{\gamma}\}\cup\{\om^m(\tau_l)\in\Omega(\delta)\}\,.
$$
\end{lemma}
\begin{proof}
Due to Lemmas~\ref{l.new} and \ref{notatzero},
\begin{equation*}
\begin{split}
\int_0^T\IP\Big(\Omega_R\cup\{I(\tau)\in Q_\gamma\}&\cup
\{W^m(\tau)\in\Omega(\delta)\}\Big)\,d\tau \\
&\le\vk_\infty(R)+\vk(\gamma^{-1};R,m)
  +\vk(\delta^{-1};R,m,\gamma)\,.
 \end{split}
\end{equation*}
Writing  the l.h.s. as $\int_0^{\nu
L}\sum_{l=0}^K\IP(\cE_l)\,dt_0$, where $\cE_l$ is defined in
terms of $\tau_l=t_0+\nu jL$, and applying the meanvalue theorem
we get the assertion.
\end{proof}

Applying the Doob inequality and Lemmas~\ref{l_polyzav},
\ref{invmeabou} to \eqref{de_ug} we get that
\begin{equation*}\begin{split}
\PP_{\Omega_R}\big(\sup_{t_l\le t \le
t_{l+1}}&|I^\nu(t)-I^\nu(t_l)|\ge P(R)\nu L+\Delta\big)\\
&\le\PP\big(\sup_{t_l\le t \le t_{l+1}}\nu\big|\int^t_{t_l}
\sigma(v(s)\,d\beta_s\big|^2\ge\Delta^2 \big) \le C_N(\nu L)^N\Delta
^{-2N}\,,
 \end{split}
\end{equation*}
 for all $N$ and $\Delta$. Choosing in this inequality $\Delta=(\nu L)^{1/3}$,
using \eqref{apr1} and  denoting
$$ Q_l=\{\sup_{t_l\le t\le
t_{l+1}}|I^\nu(t)-I^\nu(t_l)|\ge P_1(R)(\nu L)^{1/3}\}\,,
$$
where $P_1$ is a suitable polynomial, we have
\begin{equation}\label{6.n}
\PP_{\Omega_R}(Q_l)\le\vk_\infty\big((\nu L)^{-1};m\big)\,.
\end{equation}
 Let us set
$$
{\cal F}_l=\cE_l\cup Q_l\,,\quad l=0,1,\dots,K\,.
$$
Then \eqref{meas} implies the estimate
\begin{equation*}
 \begin{split}
\frac1{K}\,\sum_{l=0}^K\PP\FF_l\le
\vk_\infty(R)+\vk({\gamma}^{-1};R,m) &
+\vk(\delta^{-1};\gamma,R,m)\\&+ \vk_\infty\big((\nu
L)^{-1};m\big)=:{\kappa}\,.
\end{split}
\end{equation*}

Since $F_k(I,\varphi)$ has a polynomial growth in $I$, then
\begin{equation}
  \label{x4}
  \sum_{l=0}^K\big|(\E-\E_{\FF_l})\Upsilon^j_l\big|\le P(R)\,
  \frac1{K}\, \sum_{l=0}^K\PP\FF_l\le\kappa\quad(j=1,2,3),
\end{equation}
where we denoted by $\kappa$ another function of the same form as above.
So it remains to estimate the expectations $\E_{\FF_l}\Upsilon_l^j$ and
their sums in $l$.

First we study increments of the process $\varphi^{\nu,m}(t)$. Let
us denote
$$
\varphi^{\nu,m}(t)-\varphi^{\nu,m}(t_l)-\om^m(I^\nu(t_l))(t-t_l)=:
\Phi^\nu_l(t),
 \;\; t_l\le t\le t_{l+1}.
$$
Then
\begin{equation*}\begin{split}
\Phi^\nu_l(t)=\int_{t_l}^t\Big(\om^m(I^\nu(x))-&\om^m(I^\nu(t_l))\Big)\,dx
+ \nu\int_{t_l}^t G^m\,dx \\&+\sqrt\nu\int_{t_l}^t
g^m\,d\beta_x=:J_1+J_2+J_3\,.
\end{split}
\end{equation*}
Outside the event $\FF_l$ the term $J_1$ estimates as follows
$$
|J_1|\le P(R,m)(\nu L)^{1/3}L\,.
$$
To estimate $J_2$ and $J_3$ we assume that
\begin{equation}\label{apr2}
P(R)(\nu L)^{1/3}\le\frac12\,\gamma\,.
\end{equation}
Then outside $\FF_l$ we have
\begin{equation*}\label{x6}
|I^\nu_k(t)|\ge\frac12\gamma\quad \forall\,t\in[t_l,t_{l+1}], \ k\le
m\,,
\end{equation*}
so by Lemma~\ref{l_polyzav} and \eqref{apr2} there we have
$$
|J_2|\le\nu LC(R){\gamma}^{-1}\le C'(R)(\nu L)^{2/3}.
$$

To bound $J_3$ we introduce the stopping time
$$
t'=\min\{t\ge t_l:\, \min_{k\le
M}I_k^\nu(t)\le\gamma\;\;\text{or}\;\; |I^\nu(t)|\ge R\}\wedge
t_{l+1}\,.
$$
Then
$$
(1-\chi_{\FF_l})|J_3(t)|\le\sqrt\nu\,\Big|\int_{t_l}^{t'\wedge
t}g^m(s)\,d\beta_s\big|=:J'_3(t)\,.
$$
We have $\ \nu\E \int_{t_l}^{t'}|g^m|^2\,ds\le\nu
L\gamma^{-1}C(R,m)$. So the Doob inequality implies
\begin{equation*}\begin{split}
\ \PP_{\FF_l}\{\sup_{t_l\le t\le t_{l+1}}|J_3|\ge(\nu L)^{1/3}\}
 &\le\PP\{\sup_{t_l\le t\le t_{l+1}}
|J_3'|\ge(\nu L)^{1/3}\}\\
&\le(\nu L)^{1/3}\gamma^{-1}C(R,m).
\end{split}\end{equation*}

We have seen that
\begin{equation}
  \label{x7}
  \PP_{\FF_l}\{\Phi_l^\nu\ge P'(R,m)\nu^{1/3}L^{4/3}\}\le
   (\nu L)^{1/3}\gamma^{-1}C(R,m).
\end{equation}

Now we may estimate the terms $\Upsilon_l^j$.\smallskip

{\bf Terms $\Upsilon_j^1$.}  Since $F_k\in\Plip(h)$, then  by
\eqref{x7} `probability' $\PP_{\FF_l}$ that the integrand in
$\Upsilon_l^1$ is $\ge C(R,m)\nu^{1/3}L^{4/3}$ is bounded by $(\nu
L)^{1/3}\gamma^{-1}C(R,m)$. Since outside ${\FF_l}$ the integrand is
$\le C(R,m)$, then
$$
\sum_l\E_{\FF_l}\Upsilon_l^1\le
\nu^{1/3}C(R,m,L,\gamma).
$$
\smallskip

{\bf Terms $\Upsilon_j^2$.} By Lemma~\ref{proobraz}, outside $\FF_l$
$$
\Upsilon^2_l\le \nu\delta^{-1}C(R,m)+L\nu m^{-1}C(R).
$$
So
$$
\sum_l\E_{\FF_l}\Upsilon_l^2\le (\delta L)^{-1}C(R,m)+m^{-1}C(R).
$$\smallskip

{\bf Terms $\Upsilon_j^3$.} By Lemma~\ref{l_aver}, outside $\FF_l$
we have
 $\Upsilon_l^3\le P(R)(\nu L)^{1/3}(\nu L)$. So
$$
\sum_l\E_{\FF_l}\Upsilon_l^3\le P(R)(\nu L)^{1/3}.
$$\smallskip

Now  \eqref{x4} and the obtained estimates on the terms
$\Upsilon_l^j$ imply that
\begin{equation*}\begin{split}
\sum_l\E|\eta_l|\le\kappa+
\nu^{1/3}C(R,m,L,\gamma)+
(\delta L)^{-1}C(R,m)+m^{-1}C(R).
\end{split}
\end{equation*}
Using \eqref{6x} we arrive at the final estimate:
\begin{equation}\label{x8}
\begin{split}
{\frak A}^\nu\le\vk_\infty(R)+&C(R)m^{-1}+\nu LC(R)\\
&+    \langle\text{same terms as in the r.h.s. above}\rangle.
\end{split}
\end{equation}
It is easy to see that for any $\epsilon>0$ we can choose our
parameters in the following order
$$
R\to m\to  \gamma\to\delta\to L\to\nu\,,
$$
so that \eqref{apr1}, \eqref{apr2} hold and  the r.h.s. of
\eqref{x8} is $<\epsilon$.

Thus, we have proved

\begin{proposition}\label{p_limdrif}
The limit relation (\ref{driftlim}) holds true.
\end{proposition}

In the same way one can show that
\begin{equation}\label{momdriftlim}
{\bf E}\,\max\limits_{0\le t\le
T}\bigg|\int\limits_0^t\{F_k(I^\nu(s), \varphi^\nu(s))-\langle
F_k\rangle(I^\nu(s))\}\,ds\bigg|^4\,\to\,0\quad{\text as }\ \nu\to0\,.
\end{equation}

> From Proposition \ref{p_limdrif} taking into account the a priori
estimates we finally derive
\begin{proposition}\label{l_limdrif}
The process  \eqref{8.4}
 is a square
integrable martingale with respect to the limit measure $\QQ^0$ and
the natural filtration of $\sigma$-algebras in $C([0,\infty);h_{I+})$.
\end{proposition}
\begin{proof} Let us consider the processes
$$
N_k^{\nu_j}(\tau)=
I_k^{\nu_j}(\tau)-\int_0^\tau\langle F_k\rangle(I^{\nu_j}(s))\,ds\,,
\quad \tau\in[0,T],\;  j=1,2,\dots\,.
$$
Due to \eqref{de_ug_red} and \eqref{driftlim} we can write $N_k^{\nu_j}$ as
$$
N_k^{\nu_j}(\tau)=M_k^{\nu_j}(\tau)+\Xi_k^{\nu_j}(\tau)\,.
$$
Here $M_k^{\nu_j}=I_k^{\nu_j}-\int F_k(I^{\nu_j},\vp^{\nu_j})$ is a martingale,
 and $\Xi_k^{\nu_j}$ is a process such that
$$
\E\sup_{0\le\tau\le T}|\Xi_k^{\nu_j}(\tau)|\to0\quad\text{as}\quad\nu_j\to0\,.
$$
This convergence implies that
\begin{equation}\label{8.2}
\lim_{\nu_j\to0}{\cal L}(N_k^{\nu_j}(\cdot))=
\lim_{\nu_j\to0}{\cal L}(M_k^{\nu_j}(\cdot))
\end{equation}
in the sense that if one limit exists, then another one exists as well and the
two are equal.

Due to \eqref{6.0} and the Skorokhod theorem, we can find random processes
 $J^{\nu_j}(\tau)$ and $J(\tau)$, $0\le t\le T$, such that
${\cal L}J^{\nu_j}(\cdot)={\cal L}I^{\nu_j}(\cdot)$,  ${\cal L}J(\cdot)={\cal
Q}^0$, and
\begin{equation}\label{8.3}
J^{\nu_j}\to J\quad\text{in}\quad C([0,T],h_I)\quad\text{as}\quad \nu_j\to0\,,
\end{equation}
almost surely. By Lemma~\ref{invmeabou},
$$
\PP\{\sup_{0\le\tau\le T}|I^{\nu_j}(\tau)|\ge R\}\le CR^{-1}
$$
uniformly in $\nu_j$. Since $\langle F_k\rangle\in\Plip(h_I)$ by
 Lemmas~\ref{l_polyzav} and \ref{l_aver}, then \eqref{8.3} implies that
 the left limit in \eqref{8.2} exists and equals  \eqref{8.4}.
By Lemmas~\ref{l_polyzav} and \ref{invmeabou} the family of martingales
 $M_k^{\nu_j}(\tau)$ is uniformly integrable. Since they  converge
in distribution to the process \eqref{8.4}, then the latter is a martingale as
well.
\end{proof}

Denote $Z_k(t)\equiv I_k(t)-\int_0^t \langle F_k\rangle(I(s))ds$.
Using the same arguments as above and (\ref{momdriftlim}) we can
show that $\ Z_k(t)Z_j(t)-\int\limits_0^t \langle
A_{kj}\rangle(I(s))ds $ is a $\QQ^0$-martingale in $C([0,T);h_{I+})$.
Combining the above statement we arrive at the following theorem,
where $T>0$ and $p\ge3$ are any fixed numbers.

\begin{theorem}\label{main}
Let the process $u^\nu(t)$, $0<\nu\le1$,  be a solution of equation \eqref6 which
either is stationary in time, or satisfies the $\nu$-independent
initial condition \eqref{1.14}, where $u_0$ is non-random and
smooth. Let
 $\Psi(u^\nu(\tau))=v^\nu(\tau)=(I^\nu(\tau),\varphi^\nu(\tau))$.
Then any limiting point  $\QQ^0$  of the family ${\cal
L}\{I^\nu(\cdot)\}$ as $\nu\to0$ is a measure in $C(0,T;h^p_{I+})$ which
satisfies the estimates
$$
\int \sup_{0\le\tau\le T}|I(\tau)|^N_{h^m_I}\,\QQ_0\big(dI(\cdot))
\le C(N,m,T)<\infty\quad\forall\,N,m\in\N
$$
and solves the martingale problem in $C(0,T;h^p_I)$ with the  drift
$\langle F\rangle(I)$ and covariance $\langle A\rangle(I)$.
\end{theorem}
Let $\sigma^0(I)$ be a symmetric square root of $\langle
A\rangle(I)$ so that $(\sigma^0(I))\sigma^0(I)^t=\langle
A\rangle(I)$. We recall that $\langle A\rangle(I)$ is a positive
compact operator for each $I\in h$.

\begin{corollary}
Any limiting measure  $\QQ^0$ as in Theorem~\ref{main} is the distribution of
a solution $I(\tau)$ of the
following stochastic differential equation
\begin{equation}\label{Whi}
dI=\langle F(I)\rangle\, d\tau+\sigma^0(I)dW_\tau\,,
\end{equation}
where $W_t$ is a cylindrical Brownian motion on $h^0_I$.
\end{corollary}
\begin{proof}
Denote by $\frak h$ the Hilbert space of sequences
$\{x_1,x_2\dots,\}$ with the norm $|x|^2_{\frak
h}=\sum\limits_{j=1}^\infty j^{2(2p+4)}x^2_j$. It is easy to check
that $\frak h$ is continuously embedded in $h^p_I$, thus all the
coefficients $\langle F(I)\rangle$, $\sigma^0(I)$ and $\langle
A\rangle(I)$ are well-defined for any $I\in {\frak h}$.

By Theorem~\ref{main} and Lemma~\ref{invmeabou}
 the measure $\QQ^0$ is concentrated on
$C(0,T;h_I^{2p+4})$.
 Since this space is continuously embedded in $C(0,T;{\frak h})$, then
$\QQ^0$ is also concentrated on $C(0,T;{\frak h})$. Therefore,
$\QQ^0$ is a solution of the above limit martingale problem in the
Hilbert space
  ${\frak h}$. It
remains to use Theorem~IV.3.5 in \cite{Yor74} (also see
\cite{DZ92}).
\end{proof}
\smallskip

The limiting measure  ${\cal Q}^0$ and the
 process $I(\tau)$ inherit the uniform in $\nu$ estimates on the
 processes $I^\nu(\tau)$, obtained in Sections~\ref{s1}-\ref{s5}.
 For example,
\begin{equation}\label{remark}
\PP\{I_k(\tau)<\delta\}\le\vk(\delta^{-1};k)\;\;\text{ uniformly in
$\tau\in[0,T],\;$ for any $k\ge1$}.
\end{equation}
\noindent In particular, $\  I(\tau)\in h^p_{I+}\setminus
\partial h^p_{I+}$ a.e., for any $\tau\ge0$.\smallskip

\noindent{\bf Remark.} Equation \eqref{Whi} is the Whitham equation
for the damped-driven KdV equation \eqref{1}. Our results show that
it has a weak solution in the space $h^p_{I+}$ for a given $I(0)$
which is a deterministic vector in the space $h^\infty_{I+}=\cap
h^p_{I+}$. In fact, the same arguments apply when $I(0)$ is a random
variable in $h^P_{I+}$ such that $\E\|I(0)\|^N_{h^p_I}<\infty$,
where $N$ and $p$ are large enough.\qed\smallskip

Now we assume that $u^\nu(t)$ is a stationary solution of \eqref6.
Then the limiting process $I(\tau)$ as in \eqref{Whi} is stationary
in $\tau$. We denote $q^0={\cal L}(I(0))$ (this is a  measure on the space
$h^p_{I+}$).

\begin{theorem}\label{Haar}
Let a process $u^\nu(t)$ be a stationary solution of equations
\eqref6.
 Then

 1) for any $0\le\tau\le T$   the law of  $\varphi^\nu(\tau)$
 converges weakly as $\nu\to0$ to the Haar measure $d\varphi$
 on $\T^\infty$.

2) The law of the pair $(I^\nu(\tau),\varphi^\nu(\tau))$ converges,
  along a subsequence
$\{\nu_j\}$, corresponding to the  measure $\QQ^0$,
 to the product measure $q^0\times{d\varphi}$.

 3) For any $m$ the measure $q^{0m}={\cal L}(I^m(0))$ is absolutely
 continuous with respect to the Lebesgue measure on $\R^m_+$.
\end{theorem}

More precisely, the second assertion of the theorem means the
following. Due to \eqref{remark} the limiting measure us supported
by the Borel set $h^p\cap \{v:\,\vv_j\ne0\;\forall\,j\}$, which is
measurably isomorphic to $(h^p_{I+}\setminus \partial
h^p_{I+})\times\T^\infty$. Under this isomorphism the limiting
measure reeds as $q^0\times d\varphi$.

\begin{proof} 1) Let us fix any $m$ and take a bounded Lipschitz function
$f$,
defined on the torus $\T^m\subset\T^\infty$. Then
$$
\E f(\varphi^\nu(\tau))= \frac1{T}\int_0^T\E
f(\varphi^{\nu,m}(s))\,ds=
\frac1{T}\E\int_0^Tf(\varphi^{\nu,m}(s))\,ds,
$$
where $\varphi^{\nu,m}$ satisfies \eqref{findims}. Arguing as when
estimating  the expectation in the l.h.s. of \eqref{6x} in the proof
of Theorem~\ref{main}, we get that
$$
\E\int_0^T\big(f(\varphi^{\nu,m})-\langle f\rangle\,ds\big)
\to0\quad\text{as}\;\;\nu\to0\,.
$$
Therefore $\E f(\varphi^\nu(\tau))\to\langle f\rangle$, and the
first assertion of the theorem follows.

2) Consider an arbitrary bounded Lipschitz test function of the form
$\Phi(I,\varphi)=f(I^m)g(\varphi^m)$, $m\ge1$. We have
$$
\E\Phi(I^\nu(\tau),\varphi^\nu(\tau))=\frac{\nu}{T}\,
\E\int_0^{\nu^{-1}T}f(I^{\nu,m}(t))g(\varphi^{\nu,m}(t))\,dt.
$$
Consider  a uniform partition of the interval $(0,\nu^{-1}T)$ into
sufficiently long subintervals. As was shown in the proof of
Theorem~\ref{main}, with high probability on any subinterval of the
partition the function $I^{\nu,m}(t)$ does not deviate much from a
random constant (see \eqref{6.n}),
 while the normalised integral of
$g(\varphi^{\nu,m}(t))$ approaches the integral of $g$ against the
Haar measure (see the proof of the first assertion). Therefore when
$\nu\to0$, the r.h.s. above can be written as
$$
\Big(\frac{\nu}{T}\E \int\limits_0^{\nu^{-1}T}
f(I^{\nu,m}(s))ds\Big) \int\limits_{\mathbb
T^\infty}g(\varphi^m)d{\varphi}+o(1)=
$$
$$
=\int\limits_{h^p_I}f(I^m)\,dq^0 \int\limits_{\mathbb
T^\infty}g(\varphi^m)\,d{\varphi}+o(1)\,.
$$
  This completes the proof of 2).

  3) The vector $I^m(\tau)$ satisfies the Ito equation, given by the
  first $m$ components of \eqref{Whi}. The corresponding diffusion
  is non-degenerate by \eqref{6.cov}. Therefore by the Krylov
  theorem (see \cite{Kr77})
 for any Borel set $U\subset[\delta,\delta^{-1}]^m$,
  $\delta>0$, we have that
  \begin{equation}\label{6.est}
    q^{0m}(U)=\PP\{I^{\nu m}(t)\in U\}\le C_\delta|U|^{1/m}.
  \end{equation}
  Let us take any zero-set $Z\subset\R_+^m$ and write it as
  $$
  Z=Z_1\cup\dots\cup Z_m\cup\hat Z\,,\;\;\;\text{where}\;\;\;
  Z_j\subset\{I_j=0\}\;\;\text{and}\;\;\hat Z\subset\R^m_{>0}\,.
$$
Then $q^{0m}(Z_j)=0$ for each $j$ due to \eqref{remark}. Writing
$\hat Z=\cup_{\delta>0}Z_\delta$, where $Z_\delta=Z\cap
[\delta,\delta^{-1}]^m$, we use \eqref{6.est} to get that
$q^{0m}(\hat Z)=\lim q^{0m}(Z_\delta)=0$. So $q^{0m}(Z)=0$ and the
theorem's proof is completed.
\end{proof}

\noindent{\bf Remark.} For any $j\ge1$ the measure $q_j^0={\cal
L}(I_j(0))$ satisfies an analogy of estimate \eqref{6.est} with
$m=1$. Therefore $q^0_j=f_j(s)\,ds,\,s\ge0$, where the function
$f_j$ is bounded on segments $[\delta,\delta^{-1}]$. \medskip

\section {Proof of Lemma \ref{notatzero}}\label{proof}

 {\bf Step 1:} processes $\tilde\vv^{\gamma}_k(\tau)$.

 For $\eta_1,\eta_2\in\mathbb R^2\setminus
\{0\}$ we denote by $U(\eta_1,\eta_2)$ the element of $SO(2)$ such
that
$U(\eta_1,\eta_2)\frac{\eta_2}{|\eta_2|}=\frac{\eta_1}{|\eta_1|}$.
Note that
$U(\eta_2,\eta_1)=U(\eta_1,\eta_2)^{-1}=U(\eta_1,\eta_2)^*$.

In the fast time $\tau$  equation (\ref{3.1}) reads
\begin{equation}\label{scal3.1}
d\vv_k=\Big(\frac{1}{\nu}\,d\Psi_k(u)V(u)+A_k(v)\Big)d\tau
+\sum_jB_{kj}(v)d\beta^j_\tau\,,
\end{equation}
where we denoted
$$
A_k(v)=d\Psi_k(u){u_{xx}}+\frac{1}{2}\,\sum_{j\in\Z_0} b_j^2
d^2\Psi_k(u)[e_j,e_j]\,,\quad B_{kj}(v)=d\Psi_k(u)b_je_j\,.
$$
Let $v(\tau)=\{\vv_k(\tau),$ $k\ge1\}$ be a solution of the system
\eqref{scal3.1}${}_{k\in\N}$.

We introduce the functions
$$
\tilde A_k(\tilde \vv_k,v)=U(\tilde \vv_k,\vv_k)A_k(v),\quad \tilde
B_{kj}(\tilde \vv_k,v)=U(\tilde \vv_k,\vv_k)B_{kj}(v)\,,
$$
 smooth in $(\tilde\vv_k,\vv_k)$ from $(\mathbb
R^2\setminus\{0\})\times (\mathbb R^2\setminus\{0\})$, and consider
the additional  stochastic equation for $\tilde\vv_k(\tau)\in\R^2$:
\begin{equation}\label{bigsde}
d\tilde\vv_k=\tilde A_k(\tilde \vv_k,v)d\tau +\sum_j\tilde
B_{kj}(\tilde \vv_k,v)d\beta^j_\tau\,.
\end{equation}
Its  coefficients 
 are well defined for all non-zero $\vv_k$ and
$\tilde \vv_k$.

If $v(\tau)=\{\vv_k(\tau)$, $k\ge1\}$ is as above, then
eq.~\eqref{bigsde} with a prescribed initial data has a unique
solution, defined while
$$
|\vv_k|,\, |\tilde\vv_k|\ge c,\quad |v|_{h^1}\le C,
$$
where $c,C$ are any fixed positive constants. This solution may be
obtained as the last component of a solution $(v, \tilde\vv_k)$ of
the coupled system \eqref{scal3.1}${}_{k\in\N}$,~\eqref{bigsde}.
This system has a unique solution since \eqref{scal3.1}${}_{k\in\N}$
is equivalent to \eqref{1} (so it has a unique solution), while
\eqref{bigsde} is a Lipschitz equation on the domain, defined by the
conditions above.

For a  $\gamma\in(0,\frac12)$ we introduce the  stopping times
$\tau_i^+,\,i\ge0$ and $\tau_i^-,\,i\ge1$, where $ \tau^+_0=0\, $
and for $i\ge1$
\begin{equation*}
\begin{split}
&\tau^-_i\!=\inf\{\tau\ge \tau^+_{i-1}:\,|\vv_k(\tau)|\le\gamma\
\hbox{or }\big|v(\tau)\big|_{h^1}\ge\frac{1}{\gamma}\}\,, \\
& \tau^+_{i}=\inf\{\tau\ge \tau^-_{i}:\,|\vv_k(\tau)|\ge 2\gamma\
\hbox{and }\big|v(\tau)\big|_{h^1}\le\frac{1}{2\gamma})\}\,.
\end{split}
\end{equation*}
Note that $\tau_0^+\le\tau^-_1$, $\tau^-_i<\tau^+_i<\tau_{i+1}^-$ if
$i>0$, and $\tau_j^\pm\to\infty$ as $j\to\infty$.

Next we construct a continuous process $\tilde
\vv^{\gamma}_k(\tau)$, $\tau\ge0$. We set
$\vvv(\tau_0^+)=\vv_k(\tau_0^+)$.
 For $i=1$ we extend $\vvv(\tau)$  to the segment
 $\Delta_{i-1}:=[\tau^+_{i-1}, \tau^-_i]$ as a solution of equation
 \eqref{bigsde}, and on the segment $\Lambda_i=[\tau^-_i,\tau^+_i]$ we
 define it as  \footnote{If $\vv_k(0)=0$, then $\tau_0^+=\tau_1^-=0$ and the
   formula \eqref{neigh_zero} is not defined. But it happens with zero
   probability, and in this case we simply set  $\vvv\equiv0$.}
\begin{equation}\label{neigh_zero}
\tilde \vv^{\gamma}_k(\tau)=U(\tilde
\vv_k(\tau^-_i),\vv_k(\tau^-_i))\vv_k(\tau), \qquad\hbox{for
}\tau\in\Lambda_i\,.
\end{equation}

\begin{lemma} If $|\vvv(\tau^+_{i-1})|=|\vv_k(\tau^+_{i-1})|$ and
$\vvv$ satisfies \eqref{bigsde} on $\Delta_{i-1}$, then
  $|\vvv|=|\vv_k|$ everywhere  on that segment.
\end{lemma}
\begin{proof}
 Application of Ito's formula to the expression $\tilde
I^\gamma_k= \frac{1}{2}\, |\tilde\vv_k^{\gamma}|^2$ on the segment
$\Delta_{i-1}$  yields
$$
d\tilde I^\gamma_k=\big(\tilde \vv_k^{\gamma},\tilde A_k(\tilde
\vv_k^{\gamma},v)\big)d\tau+ \sum\limits_l \Big(\frac{1}{2}\,|\tilde
B_{kl}(\tilde \vv_k^{\gamma},v)|^2 d\tau+ (\tilde
\vv_k^{\gamma},\tilde B_{kl}(\tilde
\vv_k^{\gamma},v))d\beta_\tau^l\Big)\,.
$$
Similarly, $I_k=\frac12|\vv_k|^2$ satisfies
$$
d  I_k=\big(  \vv_k,  A_k(
 v)\big)d\tau+ \sum\limits_l\Big( \frac{1}{2}\,| B_{kl}(
 v)|^2 d\tau+ (  \vv_k, B_{kl}(
 v))d\beta_\tau^l\Big)\,.
$$
By construction, the drift and diffusion coefficients of these two
equations satisfy the relations
$$
\big(\tilde \vv_k^{\gamma},\tilde A_k(\tilde
\vv_k^{\gamma},v)\big)+\frac{1}   {2}\sum\limits_l |\tilde
B_{kl}(\tilde \vv^{\gamma}_k,v)|^2=\frac{|\tilde
\vv^{\gamma}_k|}{|\vv_k|}(\vv_k,A_k(v))+ \frac{1}{2}\sum\limits_l
|B_{kl}(v)|^2\,,
$$
$$
(\tilde \vv^{\gamma}_k,\tilde B_{kl}(\tilde
\vv^{\gamma}_k,v))=\frac{|\tilde \vv^{\gamma}_k|}{|\vv_k|}
(\vv_k,B_{kl}(v)).
$$
 For the squared  difference $(I_k-\tilde I^{\gamma}_k)^2$  we have
\begin{equation}\label{idiffsq}
\begin{array}{c}
\displaystyle d(I_k-\tilde I^{\gamma}_k)^2=\Big(2(I_k-\tilde
I^{\gamma}_k)\frac{|\vv_k|-|\tilde
\vv^{\gamma}_k|}{|\vv_k|}(\vv_k,A_k(v))+\\[2mm]
\displaystyle +\frac{(|\vv_k|-|\tilde
\vv^{\gamma}_k|)^2}{|\vv_k|^2}\sum\limits_l(\vv_k,B_{kl}(v))^2\Big)
\,d\tau+ d{\cal M}_{\tau},
\end{array}
\end{equation}
where ${\cal M}_{\tau}$ is a square integrable stochastic integral
whose structure is of no interest. Denote
$J^{\gamma}(\tau)=(I_k-\tilde I^\gamma_k)^2\big(
\big(\tau\vee\tau^+_i)\wedge\tau^-_{i+1}\big) \big)$. Since
$$
|\vv_k|-|\tilde \vv^{\gamma}_k|=2\,\frac{I_k-\tilde I^{\gamma}_k } {
|\vv_k|+|\tilde \vv^{\gamma}_k|}\,,
$$
then it follows from (\ref{idiffsq}) that $\ \E J^{\gamma}(\tau)\le
\E J^{\gamma}(0)+C(\gamma) \int\limits_0^{\tau} \E
J^{\gamma}(s)ds\,. $ As $J^\gamma(\tau^+_{i-1})=0$, then
$J^{\gamma}(\tau)\equiv0$ by the Gronwall lemma. That is,
$|\vvv|=|\vv_k|$ on $\Delta_{i-1}$.
\end{proof}

Applying this lemma with $i=1$ we see that \eqref{neigh_zero} with
$i=1$   is well defined, and
 $|\vvv|=|\vv_k|$ on $\Delta_0\cup\Lambda_1$. Repeating the
 construction above for $i=2,3,\dots$ we get a continuous process
 $\vvv(\tau), \tau\ge0$, satisfying \eqref{bigsde} on the segments $\Delta_i,
 i\ge0$, satisfying \eqref{neigh_zero} on the segments $\Lambda_i, i\ge1$,
 and such that
$$
|\vvv(\tau)|\equiv|\vv_k(\tau)|.
$$

 Let us abbreviate $U_i=U(\tilde \vv_k(\tau^-_i),\vv_k(\tau^-_i))$. Then on
the intervals $\Lambda_i$ the process $\tilde \vv^{\gamma}_k(\tau)$
satisfies the equation
$$
d\tilde\vv^{\gamma}_k(\tau)=U_i\Big(\frac{1}{\nu}d\Psi_k(u)V(v)+A_k\Big)d{\tau}
+U_iB_{kj}(v)d\beta^j_\tau.
$$
Finally, using the notation
$$
\hat A_k(\tilde \vv_k,v,t)=\left\{
\begin{array}{ll}
\tilde A_k(\tilde \vv_k,v),\quad &\tau\in
\bigcup\limits_{i}\Delta_i,\\
U_i\Big(\frac{1}{\nu}d\Psi_k(u)V(v)+A_k\Big)\quad &\tau\in
\bigcup\limits_{i}(\tau_i^-,\tau_i^+),
\end{array}
\right.
$$
and
$$
\hat B_{kj}(\tilde \vv_k,v,t)=\left\{
\begin{array}{ll}
\tilde B_{kj}(\tilde \vv_k,v),\quad &\tau\in
\bigcup\limits_{i}\Delta_i,\\
U_iB_{kj}(v)\quad &\tau\in \bigcup\limits_{i}(\tau_i^-,\tau_i^+),
\end{array}
\right.
$$
we represent $\tilde \vv^{\gamma}_k(\tau)$ as the Ito process
\begin{equation}\label{ito_repr_gamma_bis}
\tilde \vv^{\gamma}_k(\tau)=v_k(0)+\int\limits_0^{\tau}\hat
A_k(\tilde \vv_k^\gamma,v,s)ds+\int\limits_0^{\tau}\hat
B_{kj}(\tilde \vv_k^\gamma,v,s)d\beta^j_s\,.
\end{equation}
 Letting formally $\frac{|\tilde
\vv^{\gamma}_k|}{|\vv_k|}=1$ for $|\vv_k|=0$, we make the function
$\frac{|\tilde \vv^{\gamma}_k|}{|\vv_k|}\equiv1$  along all
trajectories.

By the definition of $\hat A_k$ and $\hat B_{kj}$ and by
Theorem~\ref{t2.2} the following bounds hold true with a suitable
integer $K$:
\begin{equation*}\label{aux_itoeq_bnd}
\begin{array}{ll}
|\hat A_k|\le C\big(\big|v\big|_{1}^K+1\big),\qquad
&\tau\in\bigcup\limits_i
\Delta_i\\
|\hat A_k|\le {C}\nu^{-1}\big( \big|v\big|_{1}^K+1\big),\qquad
&\tau\in\bigcup\limits_i (\tau_i^-,\tau_i^+)\\
|\hat B_k|_{h^1}\le C\big( \big|v\big|_{1}^K+1\big),\qquad
&\tau\in[0,\infty)
\end{array}
\end{equation*}
(cf. Lemma~\ref{l_polyzav}).
 Let us fix any  $\nu>0$. The family of processes
$\{\vv_k^{\gamma}(\cdot)\,,\,0<\gamma<1/2\}$ is tight in
$C(0,T;\R^2)$. This readily follows from (\ref{ito_repr_gamma_bis}),
Lemma~\ref{invmeabou} and the  estimates above.

Since $B_{kj}(v)=\beta_jd\Psi_k(u)e_j$, where $\Psi$ defines
diffeomorphisms $H^0\to h^0$ and $H^1\to h^1$, then
 the diffusion $\sum\hat B_{kj}d\beta^j$ in $\R^2$ is
non-degenerate and the corresponding diffusion matrix admits  lower
and upper bounds, uniform if $\big|v\big|_{1}\le R$ for any fixed
$R>0$.\medskip

{\bf Step 2:} Cut-of at a level $|v|_1=R$.

 Let us introduce Markov time $\bar\tau_R=\inf\{\tau\ge
0\,:\,\big|v(\tau)\big|_{1}\ge R\}$. We define the processes
$\vv_k^{R}$ equal to $\vv_k$ for $\tau\in[0,\bar\tau_R]$ and
satisfying the equation
$$
d\vv_k^{R}(\tau)=dW_\tau, \qquad \tau>\bar\tau_R\,,
$$
where $W_\tau=  \left(\begin{array}{c}
    \beta^1_\tau\\
    \beta^{-1}_\tau\\
\end{array}\right)$. Also, we define  $\tilde\vv_k^{\gamma,R}$  to be equal
  to
$\tilde\vv_k^{\gamma}$ for $\tau\in[0,\bar\tau_R]$ and for
$\tau>\bar\tau_R$ satisfying the equation
$$
d\tilde\vv_k^{\gamma,R}(\tau)=U(\vv_k(\bar\tau_R),
\tilde\vv_k^{\gamma}(\bar\tau_R))dW_{\tau}\,, \qquad
\tau>\bar\tau_R.\,.
$$
These processes have positive definite diffusion matrices uniformly
in $\gamma$ and $\nu$, and
$$
|{\tilde \vv}_k^{\gamma,R}|\equiv|\vv_k^R|\,.
$$
By Lemma \ref{invmeabou} we have
\begin{equation}\label{4.z}
\begin{split}
\PP\{\tilde\vv_k^\gamma(\tau)\not=\tilde\vv_k^{\gamma,R}(\tau)\
\text{for some } 0\le \tau\le T\}\to0,\\ \PP\{|\vv_k(\tau)|\not=|
\vv_k^R  (\tau)|\  \text{for some } 0\le \tau\le T\}\to 0
\end{split}
\end{equation}
as $R\to\infty$, uniformly in $\gamma$ and $\nu$. Therefore, it
suffices to prove the lemma for $\vv_k$ replaced by  $\vv_k^{R}$
with arbitrary $R$. \footnote{Indeed, for any $\e>0$ choosing first
$R$ so big that the
 probability in \eqref{4.z} is $<\e/2$ and choosing next $\delta=\delta(\e)$
 so small that the l.h.s. of \eqref{4.y}, evaluated for $\vv_k$ replaced by
 $\vv_k^R$,  also is $<\e/2$, we see that the l.h.s. of \eqref{4.y} is $<\e$,
 if $\delta$ is sufficiently  small. }\medskip

 {\bf Step 3:} limit $\gamma\to0$.

Denote a limiting (as $\gamma\to0$)  law of $\tilde\vv_k^{\gamma,R}$
in $C(0,T;\R^2)$ by $\tilde{\cal L}^0$, and let $\hat\vv_k(\tau)$ be
any process such that its law equals $\tilde{\cal L}^0$. By
construction, the relation holds $\tilde{\cal L}^0\{|\hat
\vv_k(\tau)|\in Q\}={\cal L}\{| \vv^R_k(\tau)|\in Q\}$ for any Borel set
$Q\subset\R$.
 So it suffices to
prove the lemma's assertion with $\vv_k$ replaced by ${\hat\vv}_k$.

 The process $\tilde
\vv^{\gamma,R}_k$ satisfies the relation
\begin{equation}\label{ito_repr_gamma_r}
\tilde \vv^{\gamma,R}_k(\tau)=\vv^R_k(0)+\int\limits_0^{\tau}\hat
A_{k,R}(\tilde \vv^{\gamma,R}_k,v,s)ds+\int\limits_0^{\tau}\hat
B_{kj,R}(\tilde \vv^{\gamma,R}_k,v,s)d\beta^j_s
\end{equation}
with
$$
\hat A_{k,R}=\left\{
\begin{array}{ll}
\hat A_{k},\quad &s\le \tau_R \\[2mm]
0,\quad &s>\tau_R
\end{array}\right.
$$
and
$$
 \hat B_{kj,R}=\left\{
\begin{array}{ll}
\hat B_{kj},\quad &s\le \tau_R\,, \\[2mm]
U\Big( \left(\begin{array}{c}
    1\\
    0\\
  \end{array}\right)\delta_{j,1}+
 \left(\begin{array}{c}
    0\\
    1\\
  \end{array}\right)\delta_{-j,1}\Big)
 &s>\tau_R\,,
\end{array}
 \right.
$$
where $U=U(\vv_k(\tau_R),\tilde\vv_k^\gamma(\tau_R))$.

Denote in (\ref{ito_repr_gamma_r}) the drift and martingale parts by
${\cal A}^\gamma(\tau)$ and ${\cal M}^\gamma(\tau)$, respectively.
Then
$$
{\cal A}^\gamma(\tau)=\int\limits_0^\tau\hat A_{k,R}(\tilde
\vv^{\gamma,R}_k,v,s)ds\,,\qquad {\cal
M}^\gamma(\tau)=\int\limits_0^\tau\hat B_{kj,R}(\tilde
\vv^{\gamma,R}_k,v,s)d\beta^j_s.
$$
Distributions of the pairs $({\cal A}^\gamma(\cdot),{\cal
M}^\gamma(\cdot))$ form a tight family of Borel measures in
$C(0,T;\mathbb R^4)$. Consider a limiting measure and represent it
as the distribution of a process
 $({\cal A}^0(\tau),{\cal M}^0(\tau))$.  Then
${\cal L}\big\{{\cal A}^0(\cdot)+{\cal M}^0(\cdot)\big\}=\tilde{\cal
L}^0$, so we can take for a process ${\hat\vv}_k$  above the process
${\hat\vv}_k(\tau)={\cal A}^0(\tau)+{\cal M}^0(\tau)$.
 Let $\tau_1$ and $\tau_2$ be arbitrary distinct point of
$[0,T]$ and $C_0$ -- any positive number.
 The set $\{\varphi\in C(0,T;\R^2)\,:\,
|\varphi(\tau_1)-\varphi(\tau_2)|\le C_0|\tau_1-\tau_2|\}$ is
closed, thus
\begin{equation}\label{ha}
\begin{split}
\limsup\limits_{\gamma\to 0}\, &\PP\{|{\cal A}^\gamma(\tau_1)-{\cal
A}^\gamma(\tau_2)|\le C_0|\tau_1-\tau_2|\}\\
&\le \PP\{|{\cal A}^0(\tau_1)-{\cal A}^0(\tau_2)|\le
C_0|\tau_1-\tau_2|\}\,.\end{split}
\end{equation}
Let us choose  $C_0=2\sup\{|A_k(v)|:\,|v|_1\le R\}$. Then
$$
|A^\gamma(\tau_1)-A^\gamma(\tau_2)|\le\frac12\, C_0|\tau_1-\tau_2|+
\nu^{-1}C(R)\big|(\cup\Lambda_j)\cap[0,T]\big|\,.
$$
Since
$$
\E\big|(\cup \Lambda_j)\cap[0,T]\big| \le \PP
\Big\{\sup_{[0,T]}|v(\tau)|_{h^1}\ge
\gamma^{-1}\Big\}+\E\int_0^T\chi_{|v_k(\tau)|\le2\gamma}\,d\tau,
$$
then it
 follows from \eqref{4.x} and
Theorem~2.2.4 in \cite{Kr77} that
 $\E |(\cup\Lambda_j)\cap[0,T]|\to
0$ as $\gamma\to 0$. Therefore the limit in the l.h.s. of \eqref{ha}
equals 1, and we conclude  that $\PP\{|{\cal A}^0(\tau_1)-{\cal
A}^0(\tau_2)|\le C_0|\tau_1-\tau_2|\}=1$. That is, ${\cal
A}^0(\tau)$ is $C_0$-Lipschitz continuous and ${\cal
A}^0(\tau)=\int\limits_0^{\tau} B^0(s)\,ds$, where $|B_0|\le C_0$.

We now turn to the martingale part. Since
$$
[0,T]\ni\tau\to {\cal M^\gamma(\tau)}\in\R^2\,,\quad 0<\gamma\le1\,,
$$
is a family of  continuous square integrable martingales with
respect to the natural filtration and uniformly bounded second
moments, then the limiting process ${\cal M}^0(\tau)$ is a
continuous square integrable martingale as well.
 Denote $\llangle{\cal
M^\gamma}\rrangle_{\tau}$ the bracket (quad\-ra\-tic
characteristics) of ${\cal M^\gamma}$. According to Corollary~VI.6.7
in \cite{Jacod-Shiryaev}, $\  \llangle{\cal
M}^0\rrangle_\tau=\lim_{\gamma\to0} \llangle{\cal
M^\gamma}\rrangle_{\tau}\,$.
 Since for $v\in\{v\,:\,|v|_{h^1}\le R \}$ it
holds
$$
c_1(\tau_1-\tau_2)|\xi|^2\le\Big(\big(\llangle{\cal
M}^\gamma\rrangle_{\tau_1}-\llangle{\cal
M}^\gamma\rrangle_{\tau_2}\big)\xi,\xi\Big)\le
c_1^{-1}(\tau_1-\tau_2)|\xi|^2\quad \forall\xi\in\R^2
$$
with some $c_1>0$, then the bracket $\llangle{\cal M}^0\rrangle$
satisfies the same estimate.
  In particular, $d\llangle{\cal M}^0\rrangle_{\tau}=a(\tau)d{\tau}$ for
  some
progressively measurable symmetric $2\times2$-matrix $a(\tau)$ such
that $c_1\,{\rm Id} \le a(\tau)\le c_1^{-1}\,{\rm Id}$, a.s. Then
$W_{\tau}=\int_0^{\tau}a^{-1/2}(s)\,d{\cal M}^0(s)$ is a Wiener
process in $\R^2$ and $ {\cal M}^0(\tau)=\int\limits_0^{\tau}
a^{1/2}(s)\,dW_s$.

 We have seen that for any $\nu>0$ and $R\ge1$
each weak limit of the family $\tilde \vv^{\gamma,R}_k(\tau)$ is an
Ito process of the form
$$
\hat \vv_k(\tau)=\hat \vv_k(0)+\int\limits_0^{\tau} B^0(s)\,ds
+\int\limits_0^{\tau} a^{1/2}(s)\,dW_s,
$$
where  $|B^0(\tau)|\le C_0$ and $c_1^{1/2}\,$Id$\le a^{1/2}(\tau)\le
c_1^{-1/2}\,$Id\  a.s., uniformly in $t$ and $\nu$. Since all the
coefficients of this equation are uniformly bounded and the
diffusion matrix  is positive definite, the desired statement
follows from Theorem~2.2.4 in \cite{Kr77}.

\section {Appendix}

Here we prove the a-priori estimates, claimed in Section~\ref{s1}.

Let $F:H^m\to\R$  be a smooth functional (for some $m\ge0$).
Applying formally  Ito's formula to $F(u(\tauu))$, where $u(\tauu)$ is
a solution, and taking the expectation we get
$$
\frac{d}{d\tauu}\,\E F(u(\tauu))=\E\langle\nabla F(u),\nu
u_{xx}-V(u)\rangle +\frac12\,\nu \sum_sb_s^2\E d^2F(u)[e_s,e_s].
$$
In particular, if $F(u)$ is an integral of motion for the KdV
equation, then $\langle\nabla F(u),V(u)\rangle=0$ and we have
\begin{equation}\label{1.1}
    \frac{d}{d\tauu}\,\E F(u(\tauu))=\nu\E\langle\nabla F(u),
u_{xx}\rangle +\frac12\,\nu \sum_sb_s^2\E d^2F(u)[e_s,e_s].
\end{equation}

Since $\|u\|_0^2$ is an integral of motion, then
$F(u)=\exp(\sigma\|u\|_0^2)$, $0<\sigma\le\frac12$, also is an
integral. We have:
$$
\nabla F(u)=2\sigma e^{\sigma\|u\|_0^2}u\,,\quad
d^2F(u)[e,e]=2\sigma e^{\sigma\|u\|_0^2}\|e\|_0^2 +4\sigma^2
e^{\sigma\|u\|_0^2}\langle u,e\rangle^2\,.
$$
So \eqref{1.1} implies that
$$
\frac{d}{d\tauu}\,\E e^{\sigma\|u\|_0^2}=-\nu\sigma\E\Big(
e^{\sigma\|u\|_0^2}\big(2\|u\|_1^2-B_0-2\sigma \sum
b_s^2u_s^2\big)\Big)\,,
$$
where for $r\ge0$ we set
\begin{equation*}
    B_r=\sum j^{2r}b_s^2\,.
\end{equation*}
Denoting $\hat B=\max b_s^2$ and choosing $\sigma\le(2\hat B)^{-1}$
we get that
$$
\frac{d}{d\tauu}\, \E e^{\sigma\|u\|_0^2}\le -\nu\sigma\E\Big(
e^{\sigma\|u\|_0^2}\big(\|u\|_1^2-B_0
\big)\Big)\le -\nu\sigma B_0\big(\E e^{\sigma\|u\|_0^2}-2e^{2\sigma
B_0}\big).
$$
So the estimate \eqref{1.4} holds for all $t\ge0$. In particular,
for each $N > 0$ we have
\begin{equation}\label{1.44}
\E \| u(\tauu) \|_0^N \le M_N = C\sigma^{-N/2} \cdot \langle \,
\text { the r.h.s. of \eqref{1.4}\,}\rangle\,.
\end{equation}

The KdV equation has infinitely many  integrals of motion $J_m
(u)$, $m\ge 0$, which can be written as
\begin{equation}\label{1.9}
J_m(u)=\|u\|_m^2+\sum_{r=3}^m\sum_\m\int C_{r,\m} u^{(m_1)}\dots
u^{(m_r)}\,dx\,.
\end{equation}
Here  the inner sum is taken over all integer $r$-vectors
$\m=(m_1,\dots m_r)$ such that $ 0\le m_j\le m-1\;\;\forall\,j\,$
and $\, m_1+\dots m_r=4+2m-2r\,$ (in particular, $J_0=\|u\|_0^2$).
E.g., see \cite{KaP}, p.~209.

Let us consider an integral as in \eqref{1.9},
$$
I=\int  u^{(m_1)}\dots u^{(m_f)}\,dx\,,\quad m_1+\dots+ m_f=M\,,
$$
where $f\ge2,\,M\ge1$ and $0\le m_j\le \mu-1$,
$\Theta:=\mu^{-1}(M+f/2-1)<2$ for some $\mu\ge2$. Then, by
H\"older's inequality,
$$
|I|\le|u^{(m_1)}|_{L_{p_1}}\dots |u^{(m_f)}|_{L_{p_f}}\,, \quad
p_j=\frac{M}{m_j}\le\infty.
$$
Applying next the Gigliardo--Nirenberg inequality we find that
\begin{equation}\label{1.09}
    |I|\le C\|u\|^\Theta_\mu\,\|u\|_0^{f-\Theta}\,.
\end{equation}
Finally, evoking the Young inequality we get that
\begin{equation}\label{1.10}
 |I|\le \delta\|u\|^2_\mu   +C_\delta
 \|u\|_0^{2\frac{f-\Theta}{2-\Theta}}\,,\quad \forall\,\delta>0\,.
\end{equation}

We have
$$
{\II}_1:=\langle\nabla J_m(u),{u_{xx}}\rangle =-2\|u\|^2_{m+1} +
\sum_{r=3}^{m+2}\sum_\m C'_\m u^{(m_1)}\dots u^{(m_r)}\,dx,
$$
where $m_1+\dots+m_r=6+2m-2r$. Due to \eqref{1.10} with
$\delta=1/2,\ f=r$ and $\mu=m+1$,
$$
{\II}_1\le-\frac32\,\|u\|^2_{m+1}+C\|u\|_0^{2\frac{r-\Theta}{2-\Theta}}
\le-\frac32\,\|u\|^2_{m+1}+C\big(1+\|u\|_0^{4(m+1)}\big).
$$

Next,
$$
d^2J_m(u)[v,v]=2\|v\|^2_m+\sum_r\sum_\m\int {C''}_\m\, v^{(m_1)}
v^{(m_2)}u^{(m_3)}\dots u^{(m_r)}\,dx\,.
$$
Hence,
$$
{\II}_2:=d^2J_m(u)[e_j,e_j]\le2j^{2m}+
|e_j|_{C^{m_1}}|e_j|_{C^{m_2}}\sum_{r,\m}\int {\hat
C}_\m|u^{(m_1)}\dots|u^{(m_{\hat r})}|\,dx\,,
$$
where $\hat r=r-2$ and $m_1+\dots+ m_{\hat r}=4+2m-2r-m_1-m_2=:\hat
M$, $\hat M\ge0$. Note that $|e_j|_{C^n}=j^n$ for each $j$ and $n$.
Assume first that $r\ge4$ and $\hat M>0$. Then \eqref{1.09} implies
that
$$
{\II}_2\le2j^{2m}+C\|u\|^\Theta_{m+1}j^{m_1+m_2}\|u\|_0^{r-2-\Theta}\,,
$$
with $\Theta=2-\frac{(3/2)r+m_1+m_2}{m+1}$.

By the Young inequality,
\begin{equation*}
\begin{split}
{\II}_2&\le2j^{2m}+\delta\|u\|^2_{m+1}+
C_\delta\big(j^{m_1+m_2}\|u\|_0^{r-2-\Theta}
\big)^{\frac2{2-\Theta}}\\
&\le2j^{2m}+\delta\|u\|^2_{m+1}+C_\delta
j^{2(m+1)}\big(1+\|u\|_0^{\frac43(m+1)}\big)\,.
\end{split}
\end{equation*}
It is easy to see that this estimate also holds for $r=4$ and for
$\hat M=0$.

Using in \eqref{1.1}  with $F=J_m$ the obtained bounds for ${\II}_1$
and ${\II}_2$ we get that
\begin{equation*}
\begin{split}
\frac{d}{d\tauu}\,\E J_m(u)&\le-\frac32\,\nu\E\|u\|^2_{m+1}+
C_1\nu\E(1+\|u\|_0^{4(m+1)})+\nu\sum|s|^{2m}b_s^2\\
+&\frac12\,\delta\nu\E\|u\|^2_{m+1}\sum b_s^2
+\frac12\,C_\delta\nu\sum
b_s^2s^{2(m+1)}\big(1+\E\|u\|_0^{\frac43(m+1)}\big)\,.
\end{split}
\end{equation*}
Choosing $\delta=B_0^{-1}$ and using \eqref{1.44} we arrive at the
estimate
$$
\frac{d}{d\tauu}\,\E J_m(u)\le-\nu\E\|u\|^2_{m+1}+\nu C_m\,,
$$
where $C_m$ depends on $B_{m+1}$ and $M_{4(m+1)}$.

Applying \eqref{1.10} with $\mu=m$ to \eqref{1.9} we see that
\begin{equation}\label{1.11}
    \frac12\,\|u\|^2_m-C(1+\|u\|_0^{4m})\le J_m(u)\le2\|u\|^2_m+
    C(1+\|u\|_0^{4m}).
\end{equation}
Therefore
$$
\frac{d}{d\tauu}\,\E J_m(u)\le-\frac{\nu}2\big(\E J_m(u)-C'_m\big)\,,
$$
where $C'_m$ depends on the same quantities as $C_m$. We get that
$$
\E J_m(u(\tauu)) \le\max(\E J_m(u(0)),  C'_m)
$$
for each $t\ge0$. Using \eqref{1.11} we obtain \eqref{1.12}.

Let us take any integers $m\ge0,\,k\ge1$. By the interpolation
inequality $\|u\|_m^k\le\|u\|_{mk}\|u\|_0^{k-1}$. Therefore
$$
\E\|u\|_m^k\le\big(\E\|u\|^2_{mk}\big)^{1/2}
\big(\E\|u\|_0^{2(k-1)}\big)^{1/2}.
$$
Using this inequality jointly with \eqref{1.44} and \eqref{1.12} we
get the estimate\eqref{1.13}.

\bibliography{meas}
\bibliographystyle{amsalpha}

\end{document}